\documentclass[12pt]{article}

\usepackage{color}
\usepackage{ulem}

\usepackage{siunitx}
\usepackage{graphicx} 
\usepackage{graphicx} 
\usepackage{csquotes}
\usepackage{ragged2e}
\usepackage{varioref}
\usepackage{float}
\usepackage[thinlines]{easytable}
\usepackage{caption}
\usepackage{csquotes}
\usepackage{amsfonts,amssymb,amsmath,exscale,relsize,hyperref,enumitem,amsthm,mathtools,physics}
\usepackage{setspace}
\usepackage{appendix}

\newtheorem{thm}{Theorem}
\newtheorem{lem}[thm]{Lemma}
\newtheorem{cor}[thm]{Corollary}
\newtheorem{rem}[thm]{Remark}

\newtheorem{solution*}{Solution}

\def\1{\mathbb{1}}

\newcommand\simiid{\stackrel{iid}{\sim}}
\newcommand\simind{\stackrel{ind}{\sim}}

\providecommand{\keywords}[1]
{
  \small	
  \textbf{\textit{Keywords---}} #1
}

\title{Asymptotic Bayes Optimality for Sparse Count Data}
\author{Sayantan Paul and Arijit Chakrabarti}
\date{}

\begin{document}

\maketitle

\begin{abstract}
We consider analyzing ``quasi-sparse" count data having an excess of near zero counts and some moderate or large counts. Data on the number of crimes or terrorist attacks in several regions of a country or several countries of a continent over a period of time, for example, conform to this general pattern.
Under the assumption that the observed counts are generated from independent Poisson distributions with unknown means, it is of interest to identify the counts corresponding to large mean parameters. For data on crimes or terrorist attacks, this would correspond to identifying the `hotspots' for such events. 
This problem can be formally stated as a problem of simultaneous testing for the means ($\theta_i$'s) of the Poisson distributions to identify the large means. For ``quasi-sparse" count data, it may be reasonable to assume that the $\theta_i$'s that are large or moderate are rather sparse in the whole mix of unknown means. The standard Bayesian approach for sparse modelling suggests use of a two-group mixture prior for the means in such a situation.
Such priors assign large probabilities for the small means and small probabilities for moderate or large means. Although intuitive, inference with such priors involves substantial computational complexity, as widely recognized in the literature. In this work we propose to use, as an alternative, a broad class of one-group global-local mixture priors for modelling the $\theta_i$'s and a multiple testing rule based on such priors for the problem at hand. Such priors are known to be computationally more amenable in high-dimensional parametric settings, and they have been used extensively in sparse normal problems. We investigate the performance of our testing rule 
from a decision theoretic point of view, when the mean parameters are truly generated from a two-group Gamma mixture prior.
We show that with respect to an additive $0-1$ loss (counting the number of misclassifications made by a multiple testing rule), the Bayes risk of our proposed multiple testing rule is asymptotically, within a constant multiple of the optimal multiple testing rule (called the Bayes Oracle rule after \cite{bogdan2011asymptotic}) in a two-group setting. The constant can be arranged to be close to $1$ for a wide range of configurations of the true data generating distribution and choices of one-group priors.
As far as we know, this is the first formal decision theoretic investigation of the use of one-group priors in the context of sparse count data.
Our work reinforces the argument for the use of appropriate one-group priors in sparse settings as an alternative to two-group mixture priors.
In our simulations, the ratio of the estimated Bayes risk of the optimal rule in the two-group setting to that of our proposed rule using one-group priors is shown to be very close to unity, for a wide range of sparsity levels.  
Our method also returns very satisfactory answers when applied in a real data example.

\keywords{Quasi-sparse count data, {Asymptotic optimality}, Bayes Oracle,
One-group prior, Gauss-Hypergeometric prior, empirical Bayes, three parameter beta normal prior, generalized double Pareto prior}
\end{abstract}

\section{Introduction}
\label{sec-1} 
In this era with enormous supply of data with a large number of variables, multiple hypothesis testing has become an area of keen interest for researchers in statistics. 
One of the key features in many high-dimensional situations, i.e., where the number of variables/parameters grow at least as fast as the sample size, is that the number of significant variables/parameters is rather small compared to the number of all possible variables/parameters. This is the so-called phenomenon of sparsity and such data appear in the context of various scientific investigations, see, e.g., \cite{buhlmann2011statistics}, \cite{mitchell1988bayesian} and  \cite{johnstone2004needles}.
The theoretical investigation of multiple testing rules in sparse high-dimensional settings has been a major focus of research over the last decade and a half, see, e.g., \cite{bogdan2008comparison}, \cite{bogdan2011asymptotic}, \cite{datta2013asymptotic}, \cite{carvalho2009handling}, \cite{carvalho2010horseshoe}, \cite{neuvial2012false} and  \cite{tang2018bayesian}. \\
\hspace*{0.5cm}
Most of these developments, however, have focused on data generated from normal distributions (e.g., observations coming from the normal means model or observations modelled by normal linear regression) or some related continuous distributions. 
However, the assumption of normality is not reasonable for many types of data and for that matter, statisticians often face problems where even the assumption of continuity of distributions become untenable.
There is a need, therefore, to develop and theoretically investigate multiple testing methods in such problems, which can be of great practical relevance as specific examples will bear out shortly. 
Our work in this paper is a modest step in that direction. \\
\hspace*{0.5cm} Suppose one is interested in modelling the counts of an event, for example, a terrorist attack, happening in several areas of a country or several countries in a continent over a certain period of time. A common feature of such data is that 
 most of the counts are close to zero (but not necessarily exactly zero)  while some are moderate or large. 
 Our interest in such data is partly motivated by the works of \cite{datta2016bayesian} and \cite{hamura2022global}, who dealt with similar settings. \cite{datta2016bayesian} referred to data of the kind just described as ``quasi-sparse" count data.
 \cite{datta2016bayesian} and \cite{hamura2022global}
 assumed that the $i$-th observation follows a Poisson distribution with mean $\theta_i$, i.e., if the observation vector is $\mathbf{Y}=(Y_1,Y_2,\cdots,Y_n)$, then $Y_i \simind Poi(\theta_i)$, where $Poi(\theta_i)$ denotes a Poisson distribution with mean $\theta_i$. A question of interest (e.g., in the context of data on terrorist attacks or any crime for that matter) would be to identify the `hotspots' i.e. areas with large $\theta_i$'s. This can be framed as a multiple testing problem for deciding which $\theta_i$'s are actually large and which are small. \cite{datta2016bayesian} and \cite{hamura2022global} followed a Bayesian viewpoint and modelled $\theta_i$'s by appropriate priors. It needs mentioning here, as pointed out in \cite{datta2016bayesian} that efforts for estimating the $\theta_i$'s corresponding to Poisson counts actually date back to the empirical Bayes approaches of \cite{robbins1956empirical} and \cite{kiefer1956consistency}, which were geared towards more general settings rather the sparse count data as in \cite{datta2016bayesian}.
 We first present a summary of the approach of \cite{datta2016bayesian}. This will be followed by a description of our work. This paper is partially based on the sumbitted thesis of the first author \cite{paulthesis2025},
 submitted to Indian Statistical Institute, Kolkata.\\
\hspace*{0.5cm}  
The most important requirement of modelling quasi-sparse count data is to be able to take into account the great abundance of near zero counts and the simultaneous occurrences of some moderate and large counts. 
In the Bayesian setting, this would amount to using some form of a mixture prior for the mean parameters.
As mentioned by \cite{datta2016bayesian}, one possible such choice is mixture of a distribution degenerate at zero with either another distribution degenerate at some non-zero value or a Gamma distribution, the mixture components arising with probabilities $1-p$ and $p$, respectively, for some small $p$.
See \cite{yang2009testing} in this context. As commented by \cite{datta2016bayesian}, such mixture modelling may be computationally unstable and also unsatisfactory in some sense in terms of performance, as illustrated in their simulations. \cite{datta2016bayesian} instead proposed to model data of such kind through a unimodal continuous prior on $\theta_i$. Specifically, they considered the following hierarchical prior: 
\begin{equation*} \label{eq:4.1.3}
	\theta_i|\lambda_i,\tau \simind{} Ga(\alpha,\lambda^2_i\tau^2), \lambda^2_i \simind{} \pi_1(\lambda^2_i), \tau^2 \sim \pi_2(\tau^2), \tag{1}
\end{equation*}
where $Ga(\alpha_1,\alpha_2)$ denotes a Gamma distribution with shape parameter $\alpha_1(>0)$ and scale parameter $\alpha_2(>0)$ with p.d.f. $f(x|\alpha_1,\alpha_2)=\frac{1}{\Gamma(\alpha_1) {\alpha}^{\alpha_1}_2}x^{\alpha_1-1} e^{-x/\alpha_2}$ for $x>0,$
while 
 $\pi_1(\cdot)$ and $\pi_2(\cdot)$ are densities for $\lambda^2_i$ and $\tau^2$, respectively. 
Here $\lambda_i$'s are known as the local shrinkage parameters and $\tau$ is termed the global shrinkage parameter. This is an example of application of an one-group global-local shrinkage prior for modelling the means of independent Poisson distributions.
In this hierarchical formulation, $\tau$ is used to induce an overall (hence global) sparsity and therefore a small value of it ensures a large mass of the distribution of $\theta_i$ around the origin. On the other hand, $\lambda_i$'s are intended to act ``locally" to modulate the shrinkage corresponding to individual parameters and hence the prior on them should have thick tails in order 
to ensure that the observations corresponding to true signals are left mostly unshrunk but the noises are indeed shrunk towards zero. See \cite{carvalho2009handling} and \cite{carvalho2010horseshoe} in this context. 
 It may be recalled that in the context of normal data, continuous one-group shrinkage priors have been used extensively.
 Some such examples of one-group global-local shrinkage priors are the three parameter beta normal prior of \cite{armagan2011generalized} (which contains the horseshoe of \cite{carvalho2009handling}, \cite{carvalho2010horseshoe}), the generalized double Pareto prior of \cite{armagan2013generalized}, and the normal-exponential-gamma prior of \cite{brown2010inference}.
 They were originally proposed as alternatives to their two-group counterparts in sparse problems under normality since such priors can substantially reduce computational complexity of doing inference using two-group priors. See, in this context, the discussions in \cite{carvalho2010horseshoe}.\\
\hspace*{0.5cm} For proving their theoretical results, \cite{datta2016bayesian} either used $\tau$ in \eqref{eq:4.1.3} as fixed or let $\tau \to 0$ as $n \to \infty$.
Integrating out $\theta_i$ and defining $\kappa_i=1/(1+\lambda^2_i\tau^2)$, the marginal distribution of $Y_i$ given $\kappa_i$ is obtained as
\begin{equation*}
	f(Y_i|\kappa_i) \propto (1-\kappa_i)^{Y_i}\kappa^{\alpha}_i,
\end{equation*}
i.e., marginally, each $Y_i$ (given $\kappa_i$) follows a negative binomial distribution with size $\alpha$ and probability of success $1-\kappa_i$. 
We note that $\theta_i|Y_i,\kappa_i \sim Ga(Y_i+\alpha,1-\kappa_i),$ and $ \mathbb{E}(\theta_i|Y_i,\kappa_i)=(1-\kappa_i)(Y_i+\alpha)$. Then $(1-\kappa_i)$ is the factor by which $(Y_i+\alpha)$ is shrunk towards zero
in the formula of the posterior mean. This is similar to the shrinkage of usual MLE towards zero in the expressions for Bayes estimates in normal models.
\cite{datta2016bayesian} assumed the Gauss Hypergeometric (GH) prior on $\kappa_i$ given by:
\begin{equation*} \label{eq:1.27}
     	\pi(\kappa_i|a_1,a_2,\tau,\gamma)=C_2 {\kappa}^{a_1-1}_i (1-\kappa_i)^{a_2-1}[1-(1-\tau^2)\kappa_i]^{-\gamma}, 0<\kappa_i<1, a_1,a_2,\tau,\gamma>0 \hspace*{0.05cm},  \tag{2}
     \end{equation*}
     where $C^{-1}_2= Beta(a_1,a_2) _2F_1(\gamma,a_1,a_1+a_2,1-\tau^2)$ is the normalizing constant,  $Beta(a_1,a_2)= \int_{0}^{1}t^{a_1-1}(1-t)^{a_2-1}dt$ is the usual beta function and $_2F_1$ the Gauss hypergeometric function (see \cite{armero1994prior} for the definition). For their theoretical analysis, \cite{datta2016bayesian} further fixed $a_1=a_2=\frac{1}{2}$.
    They established that for their proposed hierarchical form \eqref{eq:4.1.3}, with $\kappa_i$ modelled as \eqref{eq:1.27}, the marginal prior density of $\theta_i$ given $\tau$ is unbounded at origin, 
    i.e., $\theta_i$ has a \textit{pole at zero}. This ensures high probability near zero in the marginal distribution of $Y_i$. 
     They also showed that, 
     the posterior distribution of $\kappa_i$ concentrates near $0$ or $1$ depending on whether the observation is large or small, thus ensuring minimal shrinkage or heavy shrinkage for large or small observations, respectively. \\
\hspace*{0.5cm} As indicated before, in the Bayesian paradigm, the usual approach to model a vector $\boldsymbol{\theta}$ with mostly small values of coordinates is through a two-group mixture prior like a ``spike-and-slab prior".
\cite{datta2016bayesian} remarked that one such natural mixture prior for quasi-sparse Poisson counts is given by:
\begin{equation*} \label{eq:4.1.1}
	\theta_i \simiid (1-p) Ga(\alpha,\beta)+p Ga(\alpha,\beta+\delta), i=1,2,\cdots,n. \tag{3}
\end{equation*}
By keeping $\alpha$ and $\beta$ small one can ensure high concentration of a $Ga(\alpha,\beta)$ distribution near zero while letting $\delta >> \beta$ ensures that $Ga(\alpha,\beta+\delta)$ is flatter than the $Ga(\alpha,\beta)$ distribution. 
Note that, using \eqref{eq:4.1.1}, 
\begin{equation*} \label{eq:4.1.2}
	Y_i \simiid (1-p) NB(\alpha,\frac{1}{\beta+1})+p NB(\alpha,\frac{1}{\beta+\delta+1}),i=1,2,\cdots,n. \tag{4}
\end{equation*}
Here $NB(\alpha,p)$ denotes a negative binomial distribution with size $\alpha$ and probability of success $p$.
It may be noted that, by introducing a set of latent variables $\nu_i, i=1,2,\cdots,n$, such that $\nu_i=0$ indicates $H_{0i}:\theta_i \sim Ga(\alpha,\beta)$ is true and $\nu_i=1$ indicates $H_{1i}:\theta_i \sim Ga(\alpha,\beta+\delta)$ is true with $P(\nu_i=1)=p$ and $P(\nu_i=0)=1-p$, the marginal prior on $\theta_i$ becomes of the form  \eqref{eq:4.1.1}. \cite{datta2016bayesian} proposed an interesting multiple testing rule (described shortly) for testing $H_{0i}$ vs $H_{1i}$ for $i=1,2,\cdots,n$ simultaneously, using their one-group priors of the form \eqref{eq:4.1.3} and \eqref{eq:1.27} when the true $\theta_i$'s are generated from \eqref{eq:4.1.1} for Poisson count data. Before going into the details we briefly describe in the next paragraph the background and motivation behind this rule. \\
\hspace*{0.5cm} \cite{carvalho2010horseshoe}, in their seminal paper on horseshoe prior (in the context of normal means model), observed a very interesting fact. They compared the expression for the posterior mean under a spike-and-slab two-group prior (with a heavy-tailed slab part) with that under the horseshoe prior. Before explaining that,
it may be mentioned here that the horseshoe prior (to be described shortly) models the normal means as a scale mixture of normals and similarly to \eqref{eq:4.1.3}, such a mixture involves the local shrinkage parameters and global shrinkage parameter $\lambda_i$ and $\tau$, respectively. Also, under the spike-and-slab prior considered in \cite{carvalho2010horseshoe}, $\theta_i=0     $ $(H_{0i})$ with probability $(1-p)$ and $\theta_i\neq 0 $ $ (H_{1i})$ with probability $p$ and under $H_{1i}, \theta_i$ has a normal distribution with large variance and zero mean. The comparison in \cite{carvalho2010horseshoe} showed that in sparse situations (i.e. $p \approx 0$), the expected posterior shrinkage coefficient $\mathbb{E}(1-\kappa_i|\text{data})$ plays approximately the same role as the posterior inclusion probability $\mathbb{P}(\theta_i \neq 0|\text{data})$ corresponding to a two-group model. It is a standard fact that when
$\theta_i$'s are generated from a two-group model, the optimal multiple testing rule
for testing $H_{0i}$ vs $H_{1i}$ for $i=1,2,\cdots,n$ simultaneously, with respect to an additive $0-1$ loss function,
would reject $H_{0i}$ in favor of $H_{1i}$
if $\mathbb{P}(\theta_i \neq 0|\text{data}) >\frac{1}{2}$. These two facts inspired \cite{carvalho2010horseshoe} to propose a multiple testing rule based on horseshoe prior that rejects $H_{0i}$ if
$\mathbb{E}(1-\kappa_i|\mathbf{X})>\frac{1}{2}$ and apply the rule when $\theta_i$'s are truly generated from a two-group mixture prior.
They also corroborated their intuition through an extensive simulation study where the rule based on the horseshoe prior mimics the performance of the optimal Bayes rule in the two-group setting. 
\cite{datta2013asymptotic} studied this problem theoretically and proved that under additive $0-1$ loss function, the decision rule based on horseshoe prior attains the optimal Bayes risk, upto a multiplicative constant, when the data is truly generated from a two-group model.
Motivated by \cite{datta2013asymptotic}, \cite{ghosh2016asymptotic} considered a more general class of priors (originally suggested by \cite{polson2010shrink}) where the local shrinkage parameters are modelled as:
\begin{equation*} \label{eq:4.1.4}
	\pi_1(\lambda^2_i)=K (\lambda^2_i)^{-a-1} L(\lambda^2_i), \tag{5}
\end{equation*}
where $K > 0$ is the constant of proportionality and $a >0$ and $L:(0,\infty) \to (0,\infty)$ is a measurable non-constant slowly varying function in Karamata’s sense (see \cite{bingham_goldie_teugels_1987}), that is, $\frac{L(\alpha x)}{L(x)} \to 1$  as $x \to \infty$, for any $\alpha>0$.
They were able to establish that
the similar optimality results hold, in fact with better constants, for their choice of broad class of global-local priors
containing horseshoe as a special case. \\
\hspace*{0.5cm} These nice results about the optimal performance of one-group priors in two-group settings for the normal means model inspired \cite{datta2016bayesian} 
to investigate the performance of a multiple testing rule for the quasi-sparse count data (for testing $H_{0i}$ vs $H_{1i}$ simultaneously for $i=1,2,\cdots,n$) based on GH prior when the true $\theta_i$'s are generated from \eqref{eq:4.1.1}.
Further motivation for the choice of the GH prior is explained in \cite{datta2016bayesian}.
Their rule rejects $H_{0i}$ if $1-\mathbb{E}(\kappa_i|Y_i, \alpha,\tau, \gamma)>\zeta$, where $\zeta$ is a suitably chosen threshold. 
The main intuition behind this rule comes from a simple comparison of the expressions of the posterior means under the two-group model \eqref{eq:4.1.1}
and their one-group prior as done previously in \cite{carvalho2010horseshoe} in the context of normal models. This is explained in Section \ref{sec-4.3}.
In a simulation study, \cite{datta2016bayesian} showed that the performance of their proposed decision rule using GH prior with $a_1=a_2=\frac{1}{2},\gamma \neq 1$ is better than that of the horseshoe prior ($a_1=a_2=\frac{1}{2},\gamma=1$) for simultaneous testing of $H_{0i}$ vs $H_{1i}$ for $i=1,2,\cdots,n$, in terms of the misclassification probability when the data is generated from a two-group mixture. 
 They also stated an upper bound on the type I error of their decision rule when $\theta_i$'s are truly generated using \eqref{eq:4.1.1}. Their method also performed satisfactory in a real data example.
 However, they did not obtain any expression for the type II error or that of the Bayes risk based on GH prior. They also did not obtain the Bayes risk of the optimal rule (the Bayes Oracle) minimizing the Bayes risk in the two-group formulation \eqref{eq:4.1.1} based on any decision theoretic framework and for that matter, this study has not been undertaken in the literature till now 
 in the context of count data. \\
\hspace*{0.5cm}  The above discussion brings forth many interesting questions which warrant further study.
First and foremost, it would be of interest to find the expression for the Bayes risk of the optimal rule
(with respect to some natural loss functions) for multiple hypothesis testing for the $\theta_i$'s when they are indeed generated from a two-group gamma mixture as in \eqref{eq:4.1.1}. Secondly, it would be of keen interest to investigate if other choices of one-group shrinkage priors can be employed for simultaneous testing on quasi-sparse count data and if they can match (or better) the good performance/properties of inference using GH prior, as shown in this context in \cite{datta2016bayesian}. Most importantly, it will be very interesting to settle, if one-group priors can indeed be proxies of their two-group counterparts when evaluated in a decision theoretic sense and if they can attain the corresponding optimal Bayes risk, even in the context of count data of this kind. 
We have considered in this paper a broad class of global-local shrinkage priors of the form \eqref{eq:4.1.3} with the local shrinkage parameter $\lambda_i$ modelled as \eqref{eq:4.1.4}, where $a > 1$, and $K$ and $L(\cdot)$ are same as before.
This class of priors includes big subclasses of
three parameter beta normal priors 
and the generalized double Pareto prior.\\ 
\hspace*{0.5cm}We now briefly highlight 
our main contributions. 
In Section \ref{sec-4.2}, 
under the additive symmetric $0-1$ loss function,
we have obtained an asymptotic expression for the Bayes risk 
corresponding to the optimal rule (also called the Bayes oracle rule) induced by the two-group prior given in \eqref{eq:4.1.1} under some assumptions on the model parameters $p,\beta$ and $\delta$. This asymptotic framework is  somewhat similar in spirit to \cite{bogdan2011asymptotic}.
and is described in detail in Section \ref{sec-4.2}.
In Section \ref{sec-4.3},
as an alternative to the two-group approach, we propose a multiple testing rule based on our chosen class of one-group priors, to be used for simultaneous testing when the true $\theta_i$'s are generated from a two-group model. The rule is similar to that in \cite{datta2016bayesian} and is defined in \eqref{eq:4.3.6}.
Next, by letting the global shrinkage parameter $\tau$ go towards zero, 
we prove that the Bayes risk corresponding to the decision rule \eqref{eq:4.3.6} based on our proposed class of priors is asymptotically within a constant multiple of the optimal Bayes risk. The constant can be arranged to be close to $1$ for a wide range of configurations of the true data generating distribution and choices of one-group priors.
This is the first main contribution of this work and is stated in Theorem \ref{sparsecountabostun} in Subsection \ref{sec-4.4.1}. In order to obtain an expression for the Bayes risk of the decision rule \eqref{eq:4.3.6}, we needed to obtain non-trivial upper bounds for both type I and type II error probabilities of this rule. These are reported in Theorems \ref{sparsecounttypeItun} and \ref{sparsecounttypeIItun} in Subsection \ref{sec-4.4.3} respectively. Due to the unavailability of closed-form expressions of these two types of errors, we first derive some posterior concentration inequalities for $\kappa_i$ in Subsection \ref{sec-4.4.2} (see Theorems \ref{sparsecountshrinktun1}-\ref{sparsecountshrinktun3}), corresponding to our class of priors. These bounds are crucial for establishing Theorems \ref{sparsecounttypeItun} and \ref{sparsecounttypeIItun}.
Though Theorem \ref{sparsecountabostun} does not explicitly use the relation between the global shrinkage parameter $\tau$ and the level of sparsity $p$, we have observed in a detailed simulation study (see Remark \ref{sparsecountchoicetaubasedonp} in Subsection \ref{sec-4.4.1}) that the estimated misclassification probability corresponding to the decision rule \eqref{eq:4.3.6} becomes lowest when $\tau$ is of the order $p$. Hence, when 
the theoretical proportion of non-nulls is known, we recommend choosing $\tau$ to be of the order $p$. 
On the other hand, when this proportion is unknown, replacing $\tau$ by $\widehat{\tau}$ (inspired by the estimate of $p$ in \cite{yano2021minimax}) in \eqref{eq:4.3.6}, we are able to show that the Bayes risk of the empirical Bayes decision rule \eqref{eq:4.3.8} is also asymptotically within a constant multiple of the optimal Bayes risk in the two-group setting. This is formally stated in Theorem \ref{sparsecountabosEB}. 
The bounds on type I and type II error probabilities of rule \eqref{eq:4.3.8} are presented in Theorems \ref{sparsecounttypeIEB} and \ref{sparsecounttypeIIEB} in Subsection \ref{sec-4.4.4}, respectively. Theorems \ref{sparsecountabostun} and \ref{sparsecountabosEB} are the first of their kind in the literature in the context of non-normal data, as far as we know. 
\\
\hspace*{0.5cm} Although the conclusions of our main results are similar to those of Theorem 1 and Theorem 2 of \cite{ghosh2016asymptotic},
it must be emphasized here that these results are not obtained by merely following the techniques of \cite{ghosh2016asymptotic}. Non-trivial modifications of different arguments used by them and observing completely new technical facts are necessary. See in this context Remark \ref{chap4rem-4} of Section \ref{sec-4.4.2} and 
Remark \ref{rem-typeIIeb} of Section \ref{sec-4.4.4}. As far as we know, this is the first formal study of decision theoretic optimality of one-group priors in the context of sparse count data.\\
\hspace*{0.5cm} Our paper hopefully presents a strong case for application of one-group priors in this context as an alternative to two-group mixture priors.
    We should mention that since GH prior is not included in our chosen class of priors, our results do not provide any
    theoretical guarantee that the decision rule based on GH prior can also approximately mimic the optimal two-group solution in large samples.
    However, in our simulation studies, we have observed that performance of testing rules based on our prior and GH prior are quite similar to that of the optimal rule in two-group setting in terms of misclassification probability for a wide range of sparsity levels, specifically for small values of $p$. It is worthy of mention that our simulations more than the support of the theoretical guarantees of Theorem \ref{sparsecountabostun} and \ref{sparsecountabosEB} in that the ratio of estimated Bayes risks of the Oracle rule and that of our method remains very close to $1$ for a wide range of sparsity levels.
    Our method has also been applied for the analysis of a real data set and has returned pretty satisfactory results. The simulation results and the real-data application are presented in Sections \ref{chap5sim} and  \ref{chap5realdata}, respectively. 
In Section \ref{sec-4.4.5}, we present an overall discussion along with some possible extensions of this work. 
Section \ref{sec-4.4.6} contains the proofs of all theorems and other theoretical results.




\subsection{Notation}
For any two sequences of real numbers $\{a_n\}$ and $\{b_n\}$ with $b_n \neq 0$ for all $n$, $a_n \sim b_n$ implies $\lim_{n \to \infty}a_n/b_n=1$. By $a_n=O(b_n)$, and $a_n=o(b_n)$ we denote $|a_n/b_n|<M$ for all sufficiently large $n$, and $\lim_{n \to \infty}a_n/b_n=0$, respectively, $M>0$ being a global constant that is independent of $n$. Likewise, for any two positive real-valued functions $f_1(\cdot)$ and $f_2(\cdot)$ with a common domain of definition that is unbounded to the right $f_1(x) \sim f_2(x)$ denotes $\lim_{x \to \infty}f_1(x)/f_2(x)=1.$ Throughout this article, the indicator function of any set $A$ will always be denoted $\mathbf{1}\{A\}$.



\section{The Two-group prior and the Bayes Oracle}
\label{sec-4.2}
In this section we study the optimal multiple testing rule (in a decision theoretic sense) when the true $\theta_i$'s of the Poisson counts follow a two-component Gamma mixture as in \eqref{eq:4.1.1}. Our interest is in multiple testing problem of $H_{0i}:\nu_i=0$ against $H_{1i}:\nu_i=1$, for $i=1,2,\cdots,n$ as described earlier. We consider the usual $0-1$ loss function for individual tests and assume that the overall loss of a multiple testing procedure is the sum of losses corresponding to individual tests. In this way, our approach is based on an additive loss function and is similar to that of \cite{bogdan2011asymptotic} for the normal means model. Additive losses were earlier considered, e.g., by \cite{lehmann1957theory} in the context of testing. \\
\hspace*{0.5cm} Consider a multiple testing procedure for the above testing problem.
Let $t_{1i}$ and $t_{2i}$ denote the probabilities of type I and type II error of the $i^{\text{th}}$ test respectively, and are defined as
\begin{equation*}
	t_{1i}= P_{H_{0i}}(H_{0i} \hspace*{0.05cm}\text{is} \hspace*{0.05cm}\text{rejected})
\end{equation*}
and
\begin{equation*}
	t_{2i}= P_{H_{1i}}(H_{0i} \hspace*{0.05cm}\text{is} \hspace*{0.05cm}\text{accepted}).
\end{equation*}
Then the Bayes risk $R$ of the multiple testing procedure is obtained as
\begin{equation*} \label{eq:4.2.1}
	R= \sum_{i=1}^{n} [(1-p) t_{1i}+ p t_{2i}]  . \tag{6}
\end{equation*}
Recalling \cite{bogdan2011asymptotic}, 
it is easy to see that a multiple testing rule minimizing the Bayes risk $R$ is one that applies the simple Bayes classifier for each individual test. The optimal multiple testing rule, therefore, for each $i$, rejects $H_{0i}$ in favor of $H_{1i},$ if 
\begin{equation*} \label{eq:4.2.2}
	\frac{f_{H_{1i}}(Y_i)}{f_{H_{0i}}(Y_i)} > \frac{1-p}{p} \hspace*{0.05cm}, \tag{7}
\end{equation*}
where $f_{H_{1i}}$ and $f_{H_{0i}}$ are the marginal densities of $Y_i$ under the alternative and null hypotheses
respectively. Using \eqref{eq:4.1.2}, the optimal testing rule given in \eqref{eq:4.2.2} is simplified as
\begin{equation*} \label{eq:4.2.3}
	\text{reject} \hspace*{0.05cm} H_{0i} \hspace*{0.05cm} \text{if} \hspace*{0.5cm} Y_i > C \hspace*{0.05cm},  \tag{8}
\end{equation*}
where 
\begin{equation*} \label{eq:4.2.4}
	C=C_{p,\alpha,\beta,\delta}=\frac{\log (\frac{1-p}{p})+\alpha \log (\frac{\beta+\delta+1}{\beta+1})}{\log \bigg[\frac{(\beta+\delta)(\beta+1)}{\beta (\beta+\delta+1)}\bigg]} \hspace*{0.05cm}.  \tag{9}
\end{equation*}
Due to the presence of $p,\alpha,\beta$ and $\delta$, the decision rule \eqref{eq:4.2.3} is termed as \textit{Bayes Oracle}. Its performance can not be equalled by any testing rule based on finite samples if these parameters are unknown. Let us denote the Bayes risk of the Bayes Oracle as $R^{\text{BO}}_{\text{Opt}}$. \\
\hspace*{0.5cm} Note that, under the mixture model \eqref{eq:4.1.2}, the marginal distributions of $Y_i$ under both null and alternative hypotheses are independent of the choice of $i$. The threshold of the Bayes Oracle, given in \eqref{eq:4.2.4}, corresponding to a two-group prior is also independent of $i$. Therefore,
we drop $i$ from $t_{1i}$ and $t_{2i}$ and rename those as $t^{\text{BO}}_1$ and $t^{\text{BO}}_2$, respectively. 
Therefore, we can write 
\begin{equation*}\label{chap5br1}
R^{\text{BO}}_{\text{Opt}}=n[(1-p)t^{\text{BO}}_1+pt^{\text{BO}}_2].\tag{10}
\end{equation*}
As in \cite{bogdan2011asymptotic}, we want to study the behaviour of $R^{\text{BO}}_{\text{Opt}}$ under a suitable asymptotic setting. Since our basic model is non-normal, we cannot simply mimic
the setting of \cite{bogdan2011asymptotic} as in their Assumption (A).
Before writing out the asymptotic setting, we want to motivate the logic behind it. 
Note that here we want to model quasi-sparse count data with a mixture of $Ga(\alpha,\beta)$ and $Ga(\alpha,\beta+\delta)$ prior. In order to ensure that under the null the $\theta_i$'s are concentrated near zero, we let $\beta \to 0$ and keep $\alpha$ fixed at a small value. In order to accommodate possible large values under the alternative, we further assume that $\delta$ is fixed value bounded away from zero, which ensures that the prior mean and the variance are larger by a factor of magnitude under the alternative compared to the null. Finally, sparsity dictates that we assume $p \to 0$ as $n \to \infty$. Under these assumptions, calculations of Theorem \ref{sparsecountabostun} (in Subsection \ref{sec-4.4.1}) reveal that,
the error bounds corresponding to the Bayes oracle rule \eqref{eq:4.2.3} (with $C$ as in \eqref{eq:4.2.4}) are of the form, 
\begin{align*}
	t^{\text{BO}}_1 & \leq \mathbb{P}_{H_{0i}} \bigg(Y_i > \frac{\log(\frac{1}{p} ) +\alpha \log (1+\delta)}{\log (\frac{1}{\beta})} (1+o(1))\bigg) \\
	\text{and} \\
	t^{\text{BO}}_2 & = P_{H_{1i}} \bigg(Y_i \leq \frac{\log(\frac{1}{p} ) +\alpha \log (1+\delta)}{\log (\frac{1}{\beta})}(1+o(1))\bigg) \hspace{0.05cm},
\end{align*}
where the $o(1)$ term is non-random, independent of index $i$ and goes to $0$ as $n \to \infty$.
If we further assume that $\beta \propto p^{C_1}$ for some $C_1 >1$, then the type II error, $t^{\text{BO}}_2$ and eventually the power of the Bayes Oracle rule is asymptotically strictly between $0$ and $1$.
In this paper, we are interested in these situations, since otherwise,
even the Bayes Oracle rule cannot guarantee non-trivial inference.
We now state our asymptotic setting in \hyperlink{assumption1}{Assumption 1} below.
Our asymptotic analysis of the Bayes Oracle 
\eqref{eq:4.2.3} and our proposed rule is done under \hyperlink{assumption1}{Assumption 1}. This is discussed in more details later in this paper.\\
\textbf{\hypertarget{assumption1}{Assumption 1.}} $p \to 0$, $\alpha$ and $\delta$ are fixed and bounded away from both $0$ and infinity, and $\beta \propto p^{C_1}$ for some $C_1 >1$ as $n \to \infty$.\\
\hspace*{0.5cm} Under \hyperlink{assumption1}{Assumption 1}, we have, 
\begin{align*}
    t^{\text{BO}}_2 & = (\beta+\delta+1)^{-\alpha} \textrm{ and } t^{\text{BO}}_1=o(p) \textrm{ as } n \to \infty.
\end{align*}
Combining these two, we have 
\begin{align*} \label{chap5bo}
    R^{\text{BO}}_{\text{Opt}} &=np(\delta+1)^{-\alpha}(1+o_1(1)), \tag{11}
\end{align*}
 where $o_1(1)$ is non-random, independent of index $i$ and goes to $0$ as $n \to \infty$.
Detailed calculations establishing \eqref{chap5bo}
are available in the proof of Theorem \ref{sparsecountabostun} in Section \ref{sec-4.4.6}. 
In the next subsection, we motivate and describe rules based on one-group global-local priors.
\\
\section{Multiple testing rules using one-group priors}
\label{sec-4.3} 
In this section, we propose a couple of simultaneous hypothesis testing rules based on a broad class of global-local shrinkage priors to be applied in the context of quasi-sparse Poisson count data described before. We will investigate in the following sections how they perform vis-a-vis the optimal rule, namely the Bayes Oracle when the means are actually generated from a two-group prior. We first briefly recall that \cite{datta2016bayesian} had also applied a decision rule in this context using a subclass of the Gauss Hypergeometric (GH) priors described earlier. But they did not study its optimality properties in a decision theoretic framework. Our main interest lies in such a study using one-group priors in the context of count data. Towards that we consider a broad class of global-local shrinkage priors satisfying \eqref{eq:4.1.3} and \eqref{eq:4.1.4} described in Section \ref{sec-1}. Incidentally, it may be noted that a rich subclass of the Three Parameter Beta Normal (TPBN) priors (with $a_1>0$ and $a_2 > 1$) 
and Generalized Double Pareto priors (with $a_1 > 2$ and $a_2>0$) 
are contained in our general class. The forms of these two classes of priors are described 
in the Subsection \ref{formpriors}.
Note that when $\kappa_i$ given $\tau$ is modelled by GH prior as in \eqref{eq:1.27} with $\gamma=a_1+a_2$, then it becomes TPBN prior for $\kappa_i$ given $\tau$. In this sense, the GH priors stated earlier
contain the TPBN family. But \cite{datta2016bayesian} only considered a subclass of that, as mentioned before. Therefore the subclass considered by \cite{datta2016bayesian} does not intersect with the class of TPBN priors considered by us or for that matter any other class of priors considered by us. \\
\hspace*{0.5cm} We now introduce our decision rule based on our chosen class of priors. It may first be observed that 
the posterior distribution of $\theta_i$ depends only on $Y_i,\kappa_i$ and $\tau$ and  
is of the form
\begin{equation*} \label{eq:4.3.2}
	\theta_i|Y_i,\kappa_i,\tau \sim \text{Ga}(Y_i+\alpha,1-\kappa_i) \hspace*{0.05cm}.  \tag{12}
\end{equation*}
By Fubini's theorem, the posterior mean of $\theta_i$ given $Y_i$ and $\tau$ is therefore given by
\begin{equation*} \label{eq:4.3.3}
	\mathbb{E}(\theta_i|Y_i,\tau)=(1-\mathbb{E}(\kappa_i|Y_i,\tau))(Y_i+\alpha) \hspace*{0.05cm}.  \tag{13}
\end{equation*}
Under a two-group mixture prior, the posterior distribution of $\theta_i$'s given $\textbf{Y}=(Y_1,\cdots,Y_n)$ are independent and, in fact, the posterior distribution of $\theta_i$ only depends on $Y_i$. 
When $\theta_i$ is generated from the two-group prior \eqref{eq:4.1.1}, one has
\begin{equation*}  \label{eq:4.3.4}
	\mathbb{E}(\theta_i|Y_i)= [(1-w_i)\frac{\beta}{\beta+1}+w_i \frac{\beta+\delta}{\beta+\delta+1}](Y_i+\alpha) =w^{*}_i(Y_i+\alpha)\hspace*{0.05cm},  \tag{14}
\end{equation*}
where $w_i=\mathbb{P}(H_{1i}|Y_i)$ denotes the posterior probability of $H_{1i}$ being true and $w^{*}_i$ is the observation-specific shrinkage weight. Comparing \eqref{eq:4.3.3} and \eqref{eq:4.3.4}, it is apparent that the posterior shrinkage weight $1-\mathbb{E}(\kappa_i|Y_i,\tau)$ mimics the role of $w^{*}_i$ when one uses one-group prior instead of a two-group prior in a sparse situation.\\
\hspace*{0.5cm} Under a two-group framework, the optimal multiple testing rule for an additive $0-1$ loss, is given by:
\begin{equation*} 
	\text{Reject} \hspace*{0.05cm} H_{0i} \hspace*{0.05cm} \text{if} \hspace*{0.5cm} w^{}_i > \frac{1}{2}\hspace*{0.05cm}, i=1,2,\cdots,n.  
\end{equation*}
When $\beta \to 0$ and $\delta$ is fixed as $n \to \infty$, the first term in $w^{*}_i$ goes to zero while the second goes to $\frac{\delta}{\delta+1}$ and thus $w^{*}_i \to w_i\frac{\delta}{\delta+1}$. As a result, the optimal testing rule is approximately:
\begin{equation*} \label{eq:4.3.5}
	\text{Reject} \hspace*{0.05cm} H_{0i} \hspace*{0.05cm} \text{if} \hspace*{0.5cm} \frac{\delta+1}{\delta}w^{*}_i > \frac{1}{2}\hspace*{0.05cm}, i=1,2,\cdots,n.  \tag{15}
\end{equation*}
 Using this and the previous observation (connecting $1-\mathbb{E}(\kappa_i|Y_i,\tau)$  and $w^{*}_i$), a natural testing rule for simultaneous hypothesis testing using our class of priors can be taken as:
\begin{equation*} \label{eq:4.3.6}
	\text{Reject} \hspace*{0.05cm} H_{0i} \hspace*{0.05cm} \text{if} \hspace*{0.5cm} 1-\mathbb{E}(\kappa_i|Y_i,\tau) > \frac{\delta}{2(\delta+1)}\hspace*{0.05cm}, i=1,2,\cdots,n.  \tag{16}
\end{equation*}
{It may be recalled that the decision rule in \cite{datta2016bayesian} based on GH prior is similar and based on the same intuitive observation. Our decision rule is inspired by their work. \\}
\hspace*{0.5cm} Our main theoretical result (Theorem \ref{sparsecountabostun}) and the simulation study (see also Remark \ref{sparsecountchoicetaubasedonp}) show excellent performance of \eqref{eq:4.3.6} in terms of Bayes risk (with respect to additive $0-1$ loss), if, among other things, $\tau$ is chosen proportional to $p$ if $p$ is known. When $p$ is unknown, this strategy will not work. In this case, one may expect to get good results if $\tau$ is replaced by some $\widehat{\tau}$ that is of the order of the proportion of non-zero means using observed data. Keeping this in mind and inspired by the estimate $\text{max} \{1, \sum_{i=1}^{n} \mathbf{1} \{Y_i \geq 1 \}\}$ of the number of non-null effects in a Poisson count data considered by \cite{yano2021minimax} (albeit in the case of zero-inflated poisson counts), we propose to use an empirical Bayes version of the decision rule \eqref{eq:4.3.6} defined as
\begin{equation*} \label{eq:4.3.8}
    \text{Reject} \hspace*{0.05cm} H_{0i} \hspace*{0.05cm} \text{if} \hspace*{0.5cm} 1-\mathbb{E}(\kappa_i|Y_i,\widehat{\tau}) > \frac{\delta}{2(\delta+1)}\hspace*{0.05cm}, i=1,2,\cdots,n,
    \tag{17}
\end{equation*}
where $1-\mathbb{E}(\kappa_i|Y_i,\widehat{\tau})$ is the posterior shrinkage weight $1-\mathbb{E}(\kappa_i|Y_i,\tau)$ evaluated at $\tau= \widehat{\tau}$, where
\begin{equation*} \label{eq:4.3.7}
	\widehat{\tau}=\text{max} \biggl\{\frac{1}{n}, \frac{1}{n}\sum_{i=1}^{n}\mathbf{1} \{Y_i\geq {1}\} \biggr\}.  \tag{18}
\end{equation*}
The definition of $\widehat{\tau}$ ensures that
it never collapses to zero, avoiding a major concern for empirical Bayes estimates in the literature. See in this context \cite{carvalho2009handling},  \cite{scott2010bayes}, \cite{bogdan2008comparison} and  
\cite{datta2013asymptotic}. It also follows from \eqref{eq:newT-9.1} in the proof of Theorem \ref{sparsecounttypeIEB} (in Section \ref{sec-4.4.6}) that $\mathbb{E}(\frac{1}{n}\sum_{i=1}^{n}\mathbf{1} \{Y_i\geq {1}\})/p$ converges to a constant as $n \to \infty$ under \hyperlink{assumption1}{Assumption 1}. So it is expected that $\widehat{\tau}$ remains proportional to $p$ for large enough $n$ when $p \propto n^{-\epsilon}$ for $\epsilon \in (0,1)$.
\\
\hspace*{0.5cm} Before we close this section, we mention below two important technical details which will be used in the subsequent sections. \\
\hspace*{0.5cm} As mentioned before, our main interest is in the study of asymptotic optimality of rules \eqref{eq:4.3.6} and \eqref{eq:4.3.8} when applied in a two-group setting. 
Since both \eqref{eq:4.3.6} and \eqref{eq:4.3.8} depend on the posterior mean of $\kappa_i$, we first need the form of the posterior distribution of $\kappa_i$ given $Y_i$ and $\tau$.
It may be noted that using \eqref{eq:4.1.3} and \eqref{eq:4.1.4}, the posterior distribution of $\kappa_i$'s are independently distributed given $(Y_1,Y_2,\cdots,Y_n,\tau)$. The posterior distribution of $\kappa_i$ depends only on $Y_i$ given $\tau$ and has the following  form :
\begin{equation*} \label{eq:4.3.9}
	\pi(\kappa_i|Y_i,\tau) \propto \kappa^{a+\alpha-1}_i (1-\kappa_i)^{Y_i-a-1}L(\frac{1}{\tau^2}(\frac{1}{\kappa_i}-1))  \hspace*{0.1cm}, 0<\kappa_i<1 .\tag{19}
\end{equation*}
\hspace*{0.5cm} For the theoretical development of the paper, we consider slowly varying functions that satisfy \hyperlink{assumption2}{Assumption 2}
below.
\\
\textbf{\hypertarget{assumption2}{Assumption 2.}} \\
\textbf{(\hypertarget{A1}{A1})} There  exists  some $c_0(>0)$ such that $L(t) \geq c_0 \; \forall t \geq t_0$, for some $t_0>0$, which depends on both $L$ and $c_0$ with $\lim_{t \to \infty} L(t) \in (0,\infty)$. \\
\textbf{(\hypertarget{A2}{A2})} There exists some $M \in (0,\infty)$ such that $\sup_{t \in (0,\infty)} L(t) \leq M$. \\
These assumptions are quite similar to those of \cite{ghosh2017asymptotic}.
\section{Theoretical Results}
\label{sec-4.4}
 
This section is divided into several subsections, the first of which, namely Subsection \ref{sec-4.4.1}, describes the main theoretical findings of this work. These results depend on the type I and type II errors corresponding to the decision rules \eqref{eq:4.3.6} and \eqref{eq:4.3.8}. These bounds are stated in Subsections \ref{sec-4.4.3} and \ref{sec-4.4.4}, respectively. These bounds are obtained using some concentration inequalities, which are presented in Subsection \ref{sec-4.4.2}.


\subsection{Results on Asymptotic Bayes Risk under sparsity using one-group priors}
\label{sec-4.4.1}
In this subsection, we present and discuss our main optimality results (in terms of Bayes risk with respect to additive $0-1$ loss) for our proposed simultaneous testing rules when the true means are generated from a two-component mixture prior \eqref{eq:4.1.1}. Our optimality results are derived under the sparse asymptotic regime described in \hyperlink{assumption1}{Assumption 1} in Section \ref{sec-4.2}.  
Our first result is described in Theorem \ref{sparsecountabostun} below, 
proof of which is provided in Section \ref{sec-4.4.6}.
\begin{thm}
	\label{sparsecountabostun}
    Let $Y_i \sim Poi(\theta_i)$ independently for $i=1,2,\cdots,n$ and suppose each $\theta_i$ is generated from \eqref{eq:4.1.1}.
    Suppose we want to test $H_{0i}:\nu_i=0$ against $H_{1i}:\nu_i=1$, for $i=1,2,\cdots,n$ using decision rule \eqref{eq:4.3.6} induced by the class of priors \eqref{eq:4.1.3} satisfying \eqref{eq:4.1.4}, where $L(\cdot)$ satisfies  (\hyperlink{A1}{A1})
	and (\hyperlink{A2}{A2}) defined in Section \ref{sec-4.3} with $a > 1$. Also assume that $p, \beta$ and $\delta $ of the two-group model satisfy \hyperlink{assumption1}{Assumption 1}. Further assume that $\tau \to 0$ as $n \to \infty$
    . Then 
    the Bayes risk (with respect to additive $0-1$ loss) of the multiple testing rule \eqref{eq:4.3.6}, denoted $R_{\text{OG}}$, satisfies,
    \begin{align*}\label{poissonabostun1}
   1\leq  \liminf_{n \to \infty}  \frac{R_{\text{OG}}} {R^{\text{BO}}_{\text{Opt}}} \leq   \limsup_{n \to \infty}  \frac{R_{\text{OG}}} {R^{\text{BO}}_{\text{Opt}}} & \leq 
   (\delta+1)^{\alpha} \mathbb{P} \bigg(Y \leq 2a+\alpha -\frac{2(a+\alpha)}{(\delta+2)}\bigg). \tag{20}
    \end{align*}
    where $Y \sim NB(\alpha,\frac{1}{\delta+1})$.
\end{thm}
Theorem \ref{sparsecountabostun} shows that the Bayes risk ($R_{\text{OG}}$) of the decision rule \eqref{eq:4.3.6} based on one-group priors asymptotically remains within a multiplicative factor of the Bayes risk ($R^{\text{BO}}_{\text{Opt}}$) of the optimal decision rule in the two-group setting. 
In particular, we may say,
\begin{align*}
    R_{\text{OG}} =O(R^{\text{BO}}_{\text{Opt}}), \textrm{ as } n \to \infty.
\end{align*}
Although the upper bound on the limit superior of $  \frac{R_{\text{OG}}} {R^{\text{BO}}_{\text{Opt}}}$ (as in \eqref{poissonabostun1}) is larger than $1$, for reasonable choices of $a, \alpha$ and $\delta$, its values remain very close to $1$ as Table \ref{table:T-F-newratio} below shows. 
\begin{table}[h!]
	\renewcommand\thetable{1}
	\caption{Upper bound obtained in \eqref{poissonabostun1} for different choices of $a, \alpha$ and $\delta$ } 
	
	\centering 
	\begin{tabular}{c c c  c} 
		\hline
        {$a$} & {$\alpha$} & {$\delta$} & {Upper bound}  \\
		\hline 
    1.1 & 1.1 & 0.5 & 1.058 \\
    1.2 & 1.1 & 0.5 & 1.085 \\
    1.3 & 1.2 & 0.5 & 1.113 \\
    1.5 & 1.5 & 1.0 & 1.241 \\
    1.2 & 1.4 & 1.0 & 1.173 \\
    1.3 & 1.3 & 1.0 & 1.182 \\
    1.3 & 1.2 & 2.0 & 1.119 \\
    1.4 & 1.3 & 2.0 & 1.225 \\
    1.2 & 1.4 & 2.0 & 1.192 \\
	\hline
	\end{tabular}
	\label{table:T-F-newratio} 
\end{table}

\begin{rem} \label{sparsecountchoicetaubasedonp}
    The upper bound derived in Theorem \ref{sparsecountabostun} remains true for any $\tau \to 0$ as $n \to \infty$. In other words, this asymptotic result does not explicitly require using any relationship between the global shrinkage parameter $\tau$ and the level of sparsity $p$, provided they both tend to $0$ as $n \to \infty$. We may recall that in the case of modelling the means of a sparse normal means model by one-group priors, similar results of optimality can be ensured by judicious choices of global shrinkage parameter in terms of the level of sparsity, see e.g. (\cite{datta2013asymptotic}, \cite{ghosh2016asymptotic} and \cite{ghosh2017asymptotic}). Since there is a similarity in terms of the sparsity of parameters and similarity in prior structure used in our problem with those in the normal means problem, we suspect that the actual performance of our decision rule in this problem may as well depend on the choice of $\tau$ in terms of $p$.
    To investigate this for finite samples,
    we conducted an elaborate simulation study which is described below. We generate Poison count data of size $n=100, 150,$ and $200$ where the $\theta_i$'s are truly generated from a two-group model \eqref{eq:4.1.1}, with
we take $\alpha=1.5, \beta=0.005, \delta=10$. We then apply our testing rule \eqref{eq:4.3.6} on the generated data by modelling $\theta_i$'s by a member of TPBN family.
This procedure is
repeated 1000 times and the average proportion of misclassifications is reported. Five choices of $\tau$ based on $p$ are used, namely $\tau=p^2$, $\tau=\frac{p}{2}$, $\tau=p$, $\tau=2p$ and $\tau=\sqrt{p}$. The results are reported in Table \ref{table:T-F-3.20}. 
The estimated average misclassification probabilities when $\tau$ is of the order of $p$
are the smallest by a distance for all sample sizes and levels of sparsity considered in these simulations.
Based on the simulation, we recommend choosing $\tau$ to be of the order of $p$ for actual applications.
\end{rem}
\begin{table}[ht]
	\renewcommand\thetable{2}
	\caption{Average misclassification probabilities for different choices of $\tau$ based on 1000  replications } 
	
	\centering 
	\begin{tabular}{c c c c c c  c} 
		\hline
     {Sample Size} & {Sparsity level} & {$\tau=p^2$} & {$\tau=\frac{p}{2}$}& {$\tau=p$} & {$\tau= 2p$} & {$\tau=\sqrt{p}$}  \\
		\hline 
         &    0.01 & 0.092 & 0.046 & 0.044 & 0.044 & 0.138  \\
    &    0.02 & 0.097 & 0.053 & 0.052 & 0.054 & 0.157  \\
100     &   0.03 & 0.108 & 0.061 & 0.058 & 0.059 & 0.164   \\
    &    0.04 & 0.115 &0.069 & 0.067 & 0.071 & 0.172  \\
     &   0.05 & 0.121 &0.071 & 0.072 & 0.075 & 0.183   \\
\hline
        &    0.01 & 0.089 &0.043 & 0.041 &0.042 & 0.132  \\
    &    0.02 & 0.093 &0.048 & 0.049 & 0.051 & 0.151  \\
150     &   0.03 & 0.097 &0.053 & 0.052 &0.055 & 0.159   \\
    &    0.04 & 0.102 &0.064 & 0.061 &0.062 & 0.168  \\
     &   0.05 & 0.114 &0.070 & 0.068 & 0.069 & 0.175   \\
        \hline
    &    0.01 & 0.078 & 0.038 & 0.037 &0.040 & 0.124  \\
    &    0.02 & 0.081 &0.041 & 0.042 & 0.045 & 0.135  \\
200     &   0.03 & 0.085 &0.049 & 0.047 &0.048 & 0.148   \\
    &    0.04 & 0.089 &0.054 & 0.051 &0.055 & 0.161  \\
     &   0.05 & 0.093 &0.057 & 0.055 &0.057 & 0.173   \\
	\hline
	\end{tabular}
	\label{table:T-F-3.20} 
\end{table}

 In practice, $p$ would generally be not known. We suspect that if in \eqref{eq:4.3.6}, $\tau$ is replaced by a $\widehat{\tau}$ that is a reasonable estimate of $p$, then we might be able to mimic the optimality result in the previous theorem. Our next theorem confirms this intuition, where we make use of the estimator \eqref{eq:4.3.7} in replacing $\tau$ by $\widehat{\tau}$ in \eqref{eq:4.3.6}. The resulting empirical Bayesian multiple testing rule \eqref{eq:4.3.8} is shown to possess similar optimality property (in Theorem \ref{sparsecountabosEB} below). We also need a very mild additional assumption that $p \propto n^{-\epsilon}$  for some $0<\epsilon<1$ which covers most cases of theoretical and practical interest. We emphasize that $\epsilon$ need not be known for our result to go through. We only need that $\epsilon \in (0,1)$.  Proof of this result is provided in Section \ref{sec-4.4.6}. 

\begin{thm}
	\label{sparsecountabosEB}
     Let $Y_i \sim Poi(\theta_i)$ independently for $i=1,2,\cdots,n$ and suppose each $\theta_i$ is generated from \eqref{eq:4.1.1}.
    Suppose we want to test $H_{0i}:\nu_i=0$ against $H_{1i}:\nu_i=1$, for $i=1,2,\cdots,n$ using decision rule \eqref{eq:4.3.8} induced by the class of priors \eqref{eq:4.1.3} and \eqref{eq:4.1.4}, where $L(\cdot)$ satisfies  (\hyperlink{A1}{A1})
	and (\hyperlink{A2}{A2}) with $a > 1$. Also assume that $p,\beta$ and $\delta$ of the two-group model satisfy \hyperlink{assumption1}{Assumption 1} with $p \propto n^{-\epsilon}$  for some $0<\epsilon<1$. Then
    the Bayes risk (with respect to additive $0-1$ loss) of the multiple testing rule \eqref{eq:4.3.8},  denoted $R^{\text{EB}}_{\text{OG}},$ satisfies,
    \begin{align*}\label{poissonaboseb1}
     1 \leq \liminf_{n \to \infty} \frac{ R^{\text{EB}}_{\text{OG}}}{ R^{\text{BO}}_{\text{Opt}}} \leq  \limsup_{n \to \infty} \frac{ R^{\text{EB}}_{\text{OG}}}{ R^{\text{BO}}_{\text{Opt}}} & \leq 
      (\delta+1)^{\alpha} \mathbb{P} \bigg(Y \leq 2a+\alpha -\frac{2(a+\alpha)}{(\delta+2)}\bigg), \tag{21}
    \end{align*}
    where $Y \sim NB(\alpha,\frac{1}{\delta+1})$.
\end{thm}
\cite{ghosh2016asymptotic} proved a similar result for the same class of priors using an empirical Bayes estimate of $\tau$ given by  \cite{van2014horseshoe} for the normal means model. Our result confirms similar phenomenon even in the case of Poisson count data. 
\begin{rem}
    If in the definition of $\widehat{\tau}$, one uses $\frac{1}{n}\sum_{i=1}^{n}\mathbf{1} \{Y_i\geq {k}\}$ for any $k \geq 2$ instead of $k=1$, the conclusions of Theorem \ref{sparsecountabosEB} remain valid, although $K=1$ makes the algebra neater. The proof for $k \geq 2$
    is similar and hence excluded. Other choices of $k$ offer a large flexibility in the definition of $\widehat{\tau}$ in practical applications, keeping in mind that the definition/perception of the ``non-null" observations can be different in different problems. In our real data application, we have argued why such a modification would be natural and more appropriate than the definition of $\widehat{\tau}$ used in Theorem \ref{sparsecountabosEB}.  
\end{rem}

\subsection{Posterior Concentration inequalities}
\label{sec-4.4.2}
In this subsection, we describe some concentration and moment inequalities for the posterior distribution of the shrinkage coefficient $\kappa_i$'s. These are essential for obtaining upper bounds on probabilities of type I and type II errors for the decision rules \eqref{eq:4.3.6} and \eqref{eq:4.3.8}. These upper bounds are discussed in the next subsections. 
The need for finding upper bounds on the errors of both type arises since it is infeasible to find exact asymptotic orders of them. This is clear from the fact that both quantities involve the posterior distribution of $\kappa_i$ as in \eqref{eq:4.3.9} and it involves a function $L$ whose form is not explicit. All of these results require use of properties of slowly varying function  and an integral inequality, which are described and proved in Lemmas \ref{chap4lem1} and \ref{chap4lem2} of Section \ref{sec-4.4.6}. Our first result of this section follows, the proof of which is given in Section \ref{sec-4.4.6}.

\begin{thm}
	\label{sparsecountshrinktun1}
	Suppose $Y_i \sim Poi(\theta_i)$ independently for $i=1,2,\cdots,n$. Consider the one-group prior given in \eqref{eq:4.1.3} satisfying \eqref{eq:4.1.4}, where $L(\cdot)$ satisfies  (\hyperlink{A1}{A1}). Then for any fixed $\epsilon \in (0,1)$, $a>0$ and any $\tau >0$,
	\begin{equation*}
		\mathbb{P}(\kappa_i < \epsilon|Y_i, \tau)  \leq \frac{a}{c_0(a+\alpha)} (K_0^{-a}-K_1^{-a})^{-1} (\tau^2)^{a-Y_i} \bigg(\frac{\epsilon}{1-\epsilon} \bigg)^{a+\alpha} L \bigg(\frac{1}{\tau^2}\bigg) (1+K_1 \tau^2)^{(Y_i+\alpha)} (1+o(1)) \hspace{0.05cm},
	\end{equation*}
	where the $o(1)$ term is non-random, independent of index $i$ and depends only on $\tau$ such that $\lim_{\tau \to 0}o(1)=0$, and $K_0$ and $K_1$ are as in Lemma \ref{chap4lem2}.
\end{thm}
Next we describe an upper bound on the posterior mean of the shrinkage coefficient $(1-\kappa_i)$,
proof of which is provided in Section \ref{sec-4.4.6}.
\begin{thm}
	\label{sparsecountshrinktun2}
	Consider the setup of Theorem \ref{sparsecountshrinktun1} with the global-local prior of the form \eqref{eq:4.1.3} satisfying \eqref{eq:4.1.4} where $L(\cdot)$ satisfies \hyperlink{assumption2}{Assumption 2} of Section \ref{sec-4.3}. Then, for any fixed $\tau >0$ and for any $Y_i \in [0, a-1)$ and $a>1$,
	\begin{equation*}
		\mathbb{E}(1-\kappa_i|Y_i, \tau) \leq \frac{a}{c_0}(K_0^{-a}-K_1^{-a})^{-1} \tau^2 [K^{-1}+ \frac{M}{(a-Y_i-1)}] (1+K_1 \tau^2)^{(Y_i+\alpha)} \hspace{0.05cm},
	\end{equation*}
    where $K_0$ and $K_1$ are as in Lemma \ref{chap4lem2}.
\end{thm}
\begin{rem}
	\label{chap4rem-4}
	For the normal means model,  \cite{ghosh2016asymptotic} derived a result on the posterior distribution of $\kappa_i$ similar to  Theorem \ref{sparsecountshrinktun1} of this work.
    However, their result is based on the Dominated Convergence Theorem. In contrast, ours uses a lower bound on the normalizing quantity of the posterior distribution of $\kappa_i$ given $Y_i$ and $\tau$, as given in Lemma \ref{chap4lem2} and the properties of slowly varying function, as given in Lemma \ref{chap4lem1}. They also obtained a result on the posterior
    mean of $1-\kappa_i$ similar to Theorem \ref{sparsecountshrinktun2}. Under \hyperlink{assumption2}{Assumption 2} on $L(\cdot)$, they provided an upper bound on $\mathbb{E}(1-\kappa_i|Y_i, \tau)$ by using the definition, $\mathbb{E}(1-\kappa_i|Y_i, \tau)= \int_{0}^{1}  \mathbb{P}(\kappa_i < \epsilon|Y_i, \tau) d \epsilon$. The same technique can not be used in our situation as $\int_{0}^{1} (\frac{\epsilon}{1-\epsilon} )^{a+\alpha} d \epsilon$ converges only when $0< a+\alpha <1$. However, for any $a > 1$, $a+\alpha > 1$ always. As a result, we need to use a completely different argument 
    to obtain a non-trivial upper bound on the posterior expectation of the shrinkage coefficient $1-\kappa_i$ in Theorem \ref{sparsecountshrinktun2}.
\end{rem}
We now find an upper bound for $\mathbb{P}(\kappa_i> \eta|Y_i, \tau)$ in order to study the behaviour of the shrinkage coefficient for large observations. The next theorem provides a concentration inequality in this context. Proof of this result is given in Section \ref{sec-4.4.6}.
\begin{thm}
	\label{sparsecountshrinktun3}
	Under the set up of Theorem \ref{sparsecountshrinktun1}, for any fixed $\eta \in (0,1)$ and $\delta_1 \in (0,1)$ and for all sufficiently small $\tau(>0)$,
	\begin{equation*}
		\mathbb{P}(\kappa_i> \eta|Y_i, \tau) \leq \frac{(a+\alpha)}{K c_0} {\tau}^{-2a} (\frac{1-\eta}{1-\eta \delta_1})^{Y_i} (\eta \delta_1)^{-(a+\alpha)} \hspace{0.05cm} .
	\end{equation*}
\end{thm}
It may be noted that this result is true by taking any $a>0$.\\
\hspace*{0.05cm} Results obtained in Theorems \ref{sparsecountshrinktun1}-\ref{sparsecountshrinktun3} reveal some interesting characteristics of the posterior distribution of $\kappa_i$ based on our choice of priors. These results also lead to the following three corollaries, whose proofs are immediate and hence omitted. \\
\hspace*{0.05cm} 
We first have a corollary of Theorem \ref{sparsecountshrinktun1}.
\begin{cor}
	\label{cor-3}
	Under the assumptions of Theorem \ref{sparsecountshrinktun1}, $\mathbb{P}(\kappa_i \geq \epsilon|Y_i, \tau) \to 1$ as $\tau \to 0$ for any fixed $\epsilon \in (0,1)$ whenever $0\leq Y_i <a$. 
\end{cor}
\cite{datta2016bayesian} also obtained a similar result for the Gauss hypergeometric prior under quasi-sparse count data.
Thus for any fixed $0\leq Y_i <a$, for our choice of one-group global-local shrinkage priors, the posterior distribution of $\kappa_i$ concentrates near $1$ with high probability for small values of $\tau$.  \\
We next have a corollary to Theorem \ref{sparsecountshrinktun2}.
\begin{cor}
	\label{cor-4}
	Under the assumptions of Theorem \ref{sparsecountshrinktun2}, $ \mathbb{E}(1-\kappa_i|Y_i, \tau) \to 0$ as $\tau \to 0$ uniformly in $Y_i \in [0, a-1)$.
\end{cor}
This corollary ensures that in the case of global-local priors, for small values of $\tau$, the noise observations will be shrunk towards the origin.\\   
Theorem \ref{sparsecountshrinktun3} yields the following corollary.
\begin{cor}
    Under the assumptions of Theorem \ref{sparsecountshrinktun3}, $\mathbb{P}(\kappa_i > \eta|Y_i, \tau) \to 0$ as $Y_i \to \infty$ for any sufficiently small $\tau>0$.
\end{cor}
This corollary implies that for our choice of one-group global-local shrinkage priors, the posterior distribution of $\kappa_i$ concentrates near $0$ with high probability for large values of $Y_i$ if $\tau$ is small enough.
\subsection{Type I and Type II error bounds for testing the rule {\eqref{eq:4.3.6}}}
\label{sec-4.4.3}
As mentioned earlier, this subsection provides some asymptotic bounds on the probabilities of both types of errors for the decision rule \eqref{eq:4.3.6} under the assumption that $p$ is known. Theorem \ref{sparsecounttypeItun} gives a non-trivial upper bound on the probability of type I error ($t_{1i}$) for the decision rule \eqref{eq:4.3.6}. On the other hand, Theorem \ref{sparsecounttypeIItun} provides an upper bound on the probability of type II error ($t_{2i}$).
We may recall that \cite{datta2016bayesian} did not study the decision theoretic optimality of their decision rule (based on GH prior) which has the same structure as \eqref{eq:4.3.6}. They only proved that the probability of type I error goes to zero asymptotically. However, their argument has a soft spot and the upper bound for type I error also needs modification.
We feel that proof of our result (Theorem \ref{sparsecounttypeItun}) on the upper bound of type I error may be helpful in getting rid of these lacunae but we have not pursued this in this work.
Proofs of Theorems \ref{sparsecounttypeItun} and \ref{sparsecounttypeIItun} are provided in Section \ref{sec-4.4.6}.
\begin{thm}
	\label{sparsecounttypeItun}
     Let $Y_i \sim Poi(\theta_i)$ independently for $i=1,2,\cdots,n$ and suppose each $\theta_i$ is generated from \eqref{eq:4.1.1}.
    Suppose we want to test $H_{0i}:\nu_i=0$ against $H_{1i}:\nu_i=1$, for $i=1,2,\cdots,n$ using decision rule \eqref{eq:4.3.6} induced by the class of priors \eqref{eq:4.1.3} satisfying \eqref{eq:4.1.4}, where $L(\cdot)$ satisfies  (\hyperlink{A1}{A1})
	and (\hyperlink{A2}{A2}) with $a>1$. Also assume that $p, \beta,$ and $\delta $ of the two-group model satisfy \hyperlink{assumption1}{Assumption 1}. Further, assume that $\tau \to 0$ as $n \to \infty$. Then as $n \to \infty$, the probability of type I error ($t_{1i}$) of the decision rule \eqref{eq:4.3.6}, satisfies
	\begin{equation*}
		t_{1i} \equiv  t_1 \leq \frac{2 \alpha \beta}{a}.
	\end{equation*}
\end{thm}
\begin{rem}
	\label{sparsecounttypeIptun}
	Under \hyperlink{assumption1}{Assumption 1}, not only does $t_1$ go to zero as $n \to \infty$, but its rate of convergence towards zero is faster than that of $p$. This fact will be used in obtaining an asymptotic expression for the Bayes risk $R_{\text{OG}}$, reported in Theorem \ref{sparsecountabostun}.
\end{rem}
\begin{thm}
	\label{sparsecounttypeIItun}
	Consider the set-up of Theorem \ref{sparsecounttypeItun}.
    Then
    as $n \to \infty$, the probability of type II error ($t_{2i}$) of the decision rule \eqref{eq:4.3.6}, denoted $t_2$, satisfies
	\begin{equation*}
		 t_{2i} \equiv t_2\leq   
   \mathbb{P} \bigg(Y \leq 2a+\alpha -\frac{2(a+\alpha)}{(\delta+2)}\bigg) (1+o(1)),
	\end{equation*}
 where the $o(1)$ term is non-random, independent of index $i$ and depends only on $\tau$ such that $\lim_{\tau \to 0}o(1)=0$ and $Y \sim NB(\alpha,\frac{1}{\delta+1})$.
\end{thm}

\begin{rem}
\label{rem-typeIItunOG}
This is the first result of its kind in the literature on one-group priors based on quasi-sparse count data. It is important to observe that the techniques used in this result are completely different from those of Theorem 7 of \cite{ghosh2016asymptotic}, which is also related to an upper bound of the type II error for the sparse normal means model. Their result is based on the usage of the upper bound on $\mathbb{P}(\kappa_i>\eta|Y_i,\tau)$. In contrast, ours uses the definition of the type II error, followed by a careful division on the range of $Y_i$ and lastly the Dominated Convergence Theorem.
\end{rem}

\subsection{Type I and Type II error bounds corresponding to an empirical Bayes procedure}
\label{sec-4.4.4}
In this subsection, we present asymptotic bounds on the probabilities of type I and type II errors of the individual decision rules \eqref{eq:4.3.8} corresponding to the empirical Bayes approach discussed before.
These are reported in Theorems \ref{sparsecounttypeIEB} and \ref{sparsecounttypeIIEB} respectively.
Approaches involved in the proofs of these theorems are significantly different from those employed for proving Theorem \ref{sparsecounttypeItun} and Theorem \ref{sparsecounttypeIItun}. Note that, when $\tau$ is used as a tuning parameter, the decision rule \eqref{eq:4.3.6} corresponding to $i^{\text{th}}$ test depends on $i^{\text{th}}$ observation $Y_i$ only. On the other hand, 
$\widehat{\tau}$ depends on (see \eqref{eq:4.3.7})
the entire dataset $(Y_1,Y_2,\cdots,Y_n)$. Obviously, this necessitates serious changes in the approach.
Towards that, we introduce a cut-off value that is asymptotically of the order $p$. Next, we divide the range of $\widehat{\tau}$
into two mutually disjoint parts depending on whether $\widehat{\tau}$ exceeds the cut-off or not. For one part, noting that for any fixed $y \geq 0$, $\mathbb{E}(1-\kappa|y,\tau)$ is non-decreasing in $\tau$, we can use results already reported for the case when $\tau$ is used as a tuning parameter. For the other part, the form of $\widehat{\tau}$ along with the independence of $Y_i$'s come into play for obtaining asymptotic bounds on the quantities of interest. 
We now present Theorem \ref{sparsecounttypeIEB}. 

\begin{thm}
	\label{sparsecounttypeIEB}
     Let $Y_i \sim Poi(\theta_i)$ independently for $i=1,2,\cdots,n$ and suppose each $\theta_i$ is generated from \eqref{eq:4.1.1}.
    Suppose we want to test $H_{0i}:\nu_i=0$ against $H_{1i}:\nu_i=1$, for $i=1,2,\cdots,n$ using decision rule \eqref{eq:4.3.8} induced by the class of priors \eqref{eq:4.1.3} satisfying \eqref{eq:4.1.4}, where $L(\cdot)$ satisfies  (\hyperlink{A1}{A1})
	and (\hyperlink{A2}{A2}) with $a>1$. Also assume that $p, \beta $ and $\delta $ of the two-group model satisfy \hyperlink{assumption1}{Assumption 1} with $p \propto n^{-\epsilon}$, for some $0<\epsilon<1$. Then as $n \to \infty$, the probability of type I error of the decision rule \eqref{eq:4.3.8}, denoted $ t^{\text{EB}}_{1i}$, satisfies
	\begin{equation*}
		t^{\text{EB}}_{1i} \leq  \frac{2 \alpha \beta} {a}+{\alpha \beta}+ e^{-(2 \log 2-1)(1-(\beta+\delta+1)^{-\alpha})np(1+o(1))} \hspace{0.05cm},
	\end{equation*}
	where the $o(1)$ term is non-random, independent of index $i$ and tends to zero as $n \to \infty$.
\end{thm}
\begin{thm}
	\label{sparsecounttypeIIEB}
	Consider the set-up of Theorem \ref{sparsecounttypeIEB}. Then 
    as $n \to \infty$, the probability of type II error of the decision rule \eqref{eq:4.3.8}, denoted  $t^{\text{EB}}_{2i}$ satisfies
	\begin{equation*}
		   t^{\text{EB}}_{2i} \leq  
   \mathbb{P} \bigg(Y \leq 2a+\alpha -\frac{2(a+\alpha)}{(\delta+2)}\bigg) (1+o(1)),    
	\end{equation*}
	where $Y \sim NB(\alpha,\frac{1}{\delta+1})$ and
    the $o(1)$ term is non-random, independent of index $i$ and tends to zero as $n \to \infty$.
\end{thm}
\begin{rem}
\label{rem-typeIIeb}
	While proving Theorem \ref{sparsecounttypeIEB} and Theorem \ref{sparsecounttypeIIEB}, we have used some arguments similar to those of \cite{ghosh2016asymptotic}. But only these arguments are not enough to establish these results. We need to use Theorem \ref{sparsecounttypeItun} and Theorem \ref{sparsecounttypeIItun} too to complete the proofs in this subsection. Our work shows that with some non-trivial modifications, the techniques used by \cite{ghosh2016asymptotic} can be useful even if the underlying model is not normal. 
\end{rem}

\section{Simulation Results}
\label{chap5sim}
In this section, we are going to compare through simulations the performance of our decision rules \eqref{eq:4.3.6} and \eqref{eq:4.3.8}
vis-a-vis the Bayes Oracle \eqref{eq:4.2.3}, when the data is truly generated from a two-group model and the data is modelled using a member of TPBN families as the prior on the $\theta_i$'s. We will also use the Gauss Hypergeometric (GH) prior of \cite{datta2016bayesian}. 
As mentioned earlier, the class of GH prior is not included in our class of global-local priors due to their choice of hyperparameters. As a result, our theoretical results do not provide any guarantee regarding the optimality of their chosen class of priors.   
Hence, it is of interest to also study the simulation performance of their decision rule. Our comparison will be based on the average value of the proportion of misclassified hypotheses in the simulation, which is used as an estimate of the risk with respect to additive $0-1$ loss.\\
\hspace*{0.5cm} The data generating scheme is as follows. For any fixed level of sparsity $p \in (0,1)$, we generate $n=500$ independent observations $Y_1, \cdots,Y_n$, from a two-group model \eqref{eq:4.1.2}. Here we have presented our results for $p \in \{0.01,0.02,\cdots,0.2\}$.
Motivated by \cite{datta2016bayesian} and inspired by \hyperlink{assumption1}{Assumption 1} related to
Theorem \ref{sparsecountabostun}, we take $\alpha=1.3,\delta=3$ and $\beta= 0.005$ when $p \in \{0.01, \cdots, 0.09\}$ and $\beta=0.05$ when  $p \in \{0.1,\cdots,0.2\}$. All the rules under study are applied to the generated dataset of size $500$.
This procedure is repeated $1000$ times and the average proportion of misclassifications is reported. For the simulation study, we use three parameter beta normal (TPBN) prior with both hyperparameters being $1.5$ and choose the global parameter $\tau$ of the same order as $p$. We also use the same GH prior as in the study of \cite{datta2016bayesian}.
We also use the empirical Bayes estimate of $\tau$, given in \eqref{eq:4.3.7}. The results are presented in Table \ref{table:T-F-3.21}.  In the table, the column Two-group prior refers to the results derived using rule \eqref{eq:4.2.3},
TPBN-tun refers to the results obtained for rule \eqref{eq:4.3.6}, while TPBN-EB refers to those for rule \eqref{eq:4.3.8}. 
\begin{table}
	\renewcommand\thetable{3}
	\caption{Average misclassification probabilities based on 1000  replications } 
	
	\centering 
	\begin{tabular}{c c c  c  c} 
		\hline
{Sparsity level} & {Two-group prior} & {TPBN-tun} & {TPBN-EB} & {Gauss Hypergeometric}  \\
		\hline 
        0.01 & 0.022 & 0.022 & 0.023 & 0.022  \\
        0.02 & 0.027 & 0.028 & 0.029 & 0.028  \\
        0.03 & 0.030 & 0.032 & 0.034 & 0.034  \\
        0.04 & 0.035 & 0.038 & 0.039 & 0.037  \\
        0.05 & 0.042 & 0.045 & 0.048 & 0.046  \\
        0.06 & 0.055 & 0.062 & 0.065 & 0.066  \\
        0.07 & 0.060 & 0.066 & 0.069 & 0.069  \\
        0.08 & 0.070 & 0.083 & 0.088 & 0.084  \\
        0.09 & 0.082 & 0.096 & 0.101 & 0.098  \\
	0.10 & 0.088 & 0.104 & 0.110 & 0.106  \\
    0.11 & 0.103 & 0.122 & 0.129 & 0.124  \\
    0.12 & 0.106 & 0.126 & 0.134 & 0.127  \\
    0.13 & 0.108 & 0.130 & 0.139 & 0.132  \\
    0.14 & 0.115 & 0.139 & 0.147 & 0.140 \\
    0.15 & 0.118 & 0.143 & 0.151 & 0.144  \\
    0.16 & 0.119 & 0.148 & 0.155 & 0.151  \\
    0.17 & 0.123 & 0.155 & 0.161 & 0.156  \\
    0.18 & 0.125 & 0.161 & 0.166 & 0.163  \\
    0.19 & 0.128 & 0.168 & 0.174 & 0.171 \\
    0.20 & 0.134 & 0.179 & 0.185 & 0.182 \\
	\hline
	\end{tabular}
	\label{table:T-F-3.21} 
\end{table}
When the level of sparsity is small, the misclassification probabilities of decision rules based on our chosen one-group prior are quite close to those of the Oracle rule, which aligns well with the theoretical results established in this paper. The results corresponding to the GH prior are also similar to those of ours, 
which indicates that the theoretical optimality property is possibly also true for the testing rule of  \cite{datta2016bayesian}.

\section{Real data analysis}
\label{chap5realdata}
\begin{figure}
     \caption{Performance of our method on real dataset}
\centering
\includegraphics[scale=1.10]{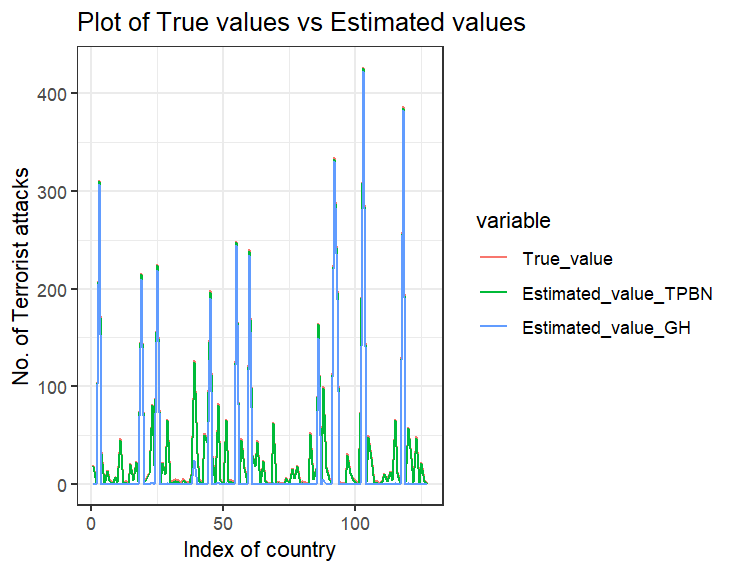}
\label{fig:1}
\end{figure}

We also apply our modelling in the Global Terrorism Dataset for the years 1970-2020, which consists of the number of terrorist attacks ($Y_i$) in different countries over the previously mentioned period. This dataset is available on the website Global Terrorism Database on request. The Figure \ref{fig:1} contains the number of terrorist attacks reported in different countries (names blinded) over the years, along with corresponding estimates using our approach and that of the GH prior-based approach. The estimates are the (empirical Bayes) posterior expected means of the rates of the terrorist attacks for different countries. We have used the empirical Bayes version of the posterior mean here as the fraction of ``nulls" is unknown here and used $\widehat{\tau}=\max \bigl\{\frac{1}{n}, \frac{1}{n}\sum_{i=1}^{n}\mathbf{1} \{Y_i\geq {100}\} \bigr\}$. The change in the definition of $\widehat{\tau}$ from \eqref{eq:4.3.8} is done due to the fact that we found it is reasonable to consider an observation coming from ``non-null" group if the average number of attacks in one year is at least $2$.
It is quite clear from the figure that 
our prior provides better estimates than those of the GH prior when the number of terrorist attacks is moderate. 
On the other hand, even if one is interested in considering the countries which are worst hit, then also, our approach provides a somewhat more accurate estimate of the true number of attacks than that of \cite{datta2016bayesian}. 
The counts for the top 10 worst hit countries along with the estimates obtained by us and that using \cite{datta2016bayesian} are illustrated in Table \ref{tab:my_label}.
 \begin{table}
   \renewcommand\thetable{4}
	\caption{Performance of our method on real dataset for worst hit countries} 
       \centering
       \begin{tabular}{c c c c}
       \hline
    {Country no.}  &  {No. of Terrorist Attacts} & {Estimate using GH}  & {Estimate using TPBN} \\
       \hline
Country 1 & 426 & 424 & 424 \\
Country 2 & 386 & 384 & 385 \\
Country 3 & 334 & 331 & 332 \\
Country 4 & 311 & 307 & 309 \\
Country 5 & 249 & 244 & 248 \\
Country 6 & 240 & 235 & 239 \\
Country 7 & 216 & 210 & 215 \\
Country 8 & 225 & 220 & 224 \\
Country 9 & 198 & 191 & 197 \\
Country 10 & 198 & 191 & 197 \\
\hline
        \end{tabular}
       \label{tab:my_label}
    \end{table}

Finally, in order to understand the role of the choice of  $\widehat{\tau}$ based on the dataset, we 
preferred simulations with different definitions of $\widehat{\tau}$. The results are presented in Table \ref{tab:chap5tab2}.
We consider countries with total number of attacks at least $100$ as belonging to non-null group.
The difference choices of $\widehat{\tau}$ use different thresholds to categorize an observation as a signal or not. We choose three values of the thresholds, namely $50, 100$ and $200$, namely, and compute the misclassification probabilities based on the decision rule \eqref{eq:4.3.8} using the corresponding $\widehat{\tau}$'s. 
It is clear from the table that the choice of $\tau$ plays a pivotal role in controlling the misclassification probability. The misclassification probability, when $\widehat{\tau}$ mimics the sample proportion of non-nulls, is significantly less compared to other choices of $\widehat{\tau}$. 
It may be recalled that a similar observation was made about the choice of $\tau$ in terms of $p$ (in case it is known) in Remark \ref{sparsecountchoicetaubasedonp}.
  \begingroup
\setlength{\tabcolsep}{10pt} 
\renewcommand{\arraystretch}{1.5} 
\begin{table}[]
        \centering
        \renewcommand\thetable{5}
	\caption{Misclassification probabilities (MP) with the choice of empirical Bayes estimate of $\tau$} 
	
	\centering
        \begin{tabular}{c c c}
        \hline
      {Non-null}  & {$\widehat{\tau}$} & {MP} \\
        \hline
        $\#\{Y_i \geq 100\}$ & $ \max\bigl\{ \frac{1}{n}, \frac{1}{n}\sum_{i=1}^{n}\mathbf{1} \{Y_i\geq {50}\} \bigr\}$ & 0.08413 \\
        $\#\{Y_i \geq 100\}$ & $ \max \bigl\{\frac{1}{n}, \frac{1}{n}\sum_{i=1}^{n}\mathbf{1} \{Y_i\geq {100}\} \bigr\}$ & 0.00744 \\
        $\#\{Y_i \geq 100\}$ & $\max \bigl\{\frac{1}{n}, \frac{1}{n}\sum_{i=1}^{n}\mathbf{1} \{Y_i\geq {200}\} \bigr\}$ & 0.02017 \\
        
        \hline    
        \end{tabular}
         \label{tab:chap5tab2}
    \end{table}

\endgroup
\section{Concluding remarks and scope for future work}
\label{sec-4.4.5}
In this paper, we study modelling of high-dimensional count data containing a lot of values near zero, namely ``quasi-sparse" count data. Assuming that observations are generated from a Poisson distribution with mean $\theta_i$ and modelling $\theta_i$ by \eqref{eq:4.1.1}, at first, we obtain an expression for the optimal Bayes risk under an additive $0-1$ loss function for simultaneous testing of $\nu_i=0$ vs $\nu_i=1$ for $i=1,2,\cdots,n$.
This was obtained under suitable assumptions on the model parameters. Next, our focus is to study the optimality property of a decision rule based on a class of one-group priors when the true data is generated from a two-group model. In this context, we have been able to establish that irrespective of whether the level of sparsity is known or not, our decision rules using global-local priors are asymptotically within a multiplicative constant of the optimal rule in the two-group setting in terms of
Bayes risk.
To the best of our knowledge, these results are the first of their kind in the literature on sparse count data.\\
\hspace*{0.5cm} Throughout our calculations, we either treat $\tau$ as a tuning parameter or use an empirical Bayes version of it. A question that arises naturally is whether the same optimality property holds for a similar decision rule when one uses a non-degenerate prior on $\tau$. Questions can also be asked regarding the choice of the prior on $\tau$ or the range on which it has to be defined.
Though we have not discussed these here, we hope techniques used in \cite{paul2025posterior} might be helpful for this purpose also.\\
\hspace*{0.5cm} \cite{hamura2022global} also studied a class of global-local priors for modeling quasi-sparse count data. They proposed two priors, namely, inverse gamma (IG) prior and extremely heavy-tailed (EH) prior, and proved that both of them possess the \textit{tail robustness} property introduced by  \cite{carvalho2010horseshoe}. However, they kept the global parameter fixed and
did not study the asymptotic optimality in terms of the Bayes risk of the decision rule induced by their chosen class of priors. We expect that by letting the global parameter go to zero as the sample size increases and using arguments similar to ours, one may be able to establish results similar to Theorem \ref{sparsecountabostun} and Theorem \ref{sparsecountabosEB} of this work. We also expect that our techniques will be handy in the case of the Gauss Hypergeometric prior used by \cite{datta2016bayesian} for providing a non-trivial bound on the probability of type II error and the corresponding Bayes risk of their decision rule, irrespective of the underlying sparsity to be known or unknown. It may be recalled that the simulation performance of the GH prior is very encouraging in terms of Bayes risk. \\
\hspace*{0.5cm} As mentioned by \cite{hamura2022global},
one possible drawback of the Poisson modelling is that for the Poisson distribution, the mean and the variance are the same, a situation which does not always hold, specially in sparse count data. One possible remedy is to use the Negative Binomial distribution to model such data and try to provide answer to the same questions discussed here. However, instead of considering a particular distribution, one may consider a general class of distributions and study such optimality for that class. We have already obtained some results in this context when observations are generated from a subclass of one parameter exponential family distributions containing Normal, Poisson and negative binomial as special cases. Interestingly, these three also belong to the  \textit{Natural Exponential Family with Quadratic Variance function}, or in short, NEF-QVF, proposed by \cite{morris1982natural}. However, optimality results in terms of Bayes risk are still not established for the class of \cite{morris1982natural}. This will be studied elsewhere.

\section{Proofs}
\label{sec-4.4.6}
\subsection{Proofs of Theorems}
Before proving the main results of this paper, first we state and prove two lemmas related to the behaviour of slowly varying functions $L$. These are very useful for the proof of the theorems that follow.
\begin{lem}
	\label{chap4lem1}
	Suppose $L$ is a slowly varying function L. Then\vspace{2mm}
	
	\begin{enumerate}
		\item \label{chap4L-1.1}  $L^{\alpha}$ is slowly varying for all $\alpha \in \mathbb{R}$.\vspace{2mm}
		
		\item \label{chap4L-1.2} $\frac{\log L(x)}{\log x} \to 0$ as $x\to \infty$.\vspace{2mm}
		
		\item \label{chap4L-1.3} for every $\alpha >0, \ x^{-\alpha} L(x) \to 0$ and $x^{\alpha}L(x) \to \infty$ as $x\to \infty$.\vspace{2mm}
		
		\item \label{chap4L-1.4} for $ \alpha<-1,-\frac{\int_x^{\infty}t^{\alpha}L(t)dt}{x^{\alpha+1}L(x)} \to \frac{1}{\alpha+1}$ as $x\to \infty$.\vspace{2mm}
		
		\item \label{chap4L-1.5} there exists a global constant $ A_0>0$ such that, for any $\alpha>-1,$ $\frac{\int_{A_0}^{x}t^{\alpha}L(t)dt}{x^{\alpha+1}L(x)} \to \frac{1}{\alpha+1}$ as $x\to \infty.$	
	\end{enumerate}

\end{lem}

\begin{proof}
    See \cite{bingham_goldie_teugels_1987}.
\end{proof}

\begin{lem}
	\label{chap4lem2}
	Let $L:(0,\infty) \to (0,\infty) $ be a measurable function satisfying \hyperlink{assumption2}{Assumption 2} of Section \ref{sec-4.3} and $a$ be any positive real number. Suppose $\tau \to 0$ as $n \to \infty$. Then for any $y =0,1,2,\cdots$, $\alpha>0$, and fixed $K_1 >K_0>\max\{1,t_0\}$,
 
	\begin{equation*}
		\int_0^{1}	u^{a+\alpha-1}(1-u)^{y-a-1}L\bigg(\frac{1}{\tau^2}(\frac{1}{u}-1)\bigg)du \geq \frac{c_0(K_0^{-a}-K_1^{-a})}{a}
		(\tau^2)^{y-a}(1+K_1\tau^2)^{-(y+\alpha)} \hspace{0.05cm},
	\end{equation*}
    where $t_0$ is as in \hyperlink{assumption2}{Assumption 2} of Section \ref{sec-4.3}.
\end{lem}
\begin{proof}
	Let $I= \int_{0}^{1}	u^{a+\alpha-1}(1-u)^{y-a-1}L\bigg(\frac{1}{\tau^2}(\frac{1}{u}-1)\bigg)du$. Then with the change of variable $t=\frac{1}{\tau^2}(\frac{1}{u}-1)$, we have
	\begin{equation*} \label{chap4eq:L-2.1}
		I= (\tau^2)^{y-a} \int_{0}^{\infty} (1+t \tau^2)^{-(y+\alpha)} t^{y-a-1} L(t) dt \hspace{0.05cm}. \tag{22}
	\end{equation*}
	By \hyperlink{assumption2}{Assumption 2} on $L(\cdot)$, $\exists K_0 (>1)$ such that $L(t) \geq c_0$ if $t \geq K_0 \geq t_0$. Hence,
    we obtain
	\begin{align*} \label{chap4eq:L-2.2}
		\int_{0}^{\infty} (1+t \tau^2)^{-(y+\alpha)} t^{y-a-1} L(t) dt &\geq 
        \int_{K_0}^{\infty} (1+t \tau^2)^{-(y+\alpha)} t^{y-a-1} L(t) dt \\
        &  \geq c_0 \int_{K_0}^{\infty} (1+t \tau^2)^{-(y+\alpha)} t^{-a-1} dt \\
        &\geq c_0
        \int_{K_0}^{K_1} (1+t \tau^2)^{-(y+\alpha)} t^{-a-1} dt \\
        & \geq \frac{c_0(K_0^{-a}-K_1^{-a})}{a} (1+K_1 \tau^2)^{-(y+\alpha)}. \tag{23}
	\end{align*}
	Finally, the result follows from \eqref{chap4eq:L-2.1} and \eqref{chap4eq:L-2.2}.
\end{proof}
\begin{proof}[Proof of Theorem \ref{sparsecountshrinktun1}]
	Fix any $\epsilon \in (0,1)$. Then
	\begin{align*}
		\mathbb{P}(\kappa_i < \epsilon|Y_i, \tau) &= \frac{\int_{0}^{\epsilon}	{\kappa_i^{a+\alpha-1}}(1-\kappa_i)^{Y_i-a-1}L\bigg(\frac{1}{\tau^2}(\frac{1}{\kappa_i}-1)\bigg)d\kappa_i}{\int_{0}^{1}	\kappa_i^{a+\alpha-1}(1-\kappa_i)^{Y_i-a-1}L\bigg(\frac{1}{\tau^2}(\frac{1}{\kappa_i}-1)\bigg)d\kappa_i} \hspace{0.05cm}.
	\end{align*}
	Now using the change of variable $t=\frac{1}{\tau^2}(\frac{1}{\kappa_i}-1)$ to the numerator and denominator of the previous quantity and applying \eqref{chap4eq:L-2.2} of Lemma \ref{chap4lem2} for the denominator, we have
	\begin{equation*} \label{eq:T-3.1}
		\mathbb{P}(\kappa_i < \epsilon|Y_i, \tau)  \leq \frac{a}{c_0} (K_0^{-a}-K_1^{-a})^{-1} (1+K_1 \tau^2)^{(Y_i+\alpha)} \int_{\frac{1}{\tau^2}(\frac{1}{\epsilon}-1)} ^{\infty} (1+t \tau^2)^{-(Y_i+\alpha)} t^{Y_i-a-1} L(t) dt \hspace{0.05cm} . \tag{24}
	\end{equation*}
	 Also, note that
\begin{align*}\label{eq:T-3.2}
    \int_{\frac{1}{\tau^2}(\frac{1}{\epsilon}-1)} ^{\infty} (1+t \tau^2)^{-(Y_i+\alpha)} t^{Y_i-a-1} L(t) dt&=
   (\tau^2)^{-Y_i-\alpha} \int_{\frac{1}{\tau^2}(\frac{1}{\epsilon}-1)} ^{\infty} \bigg(\frac{t \tau^2}{1+t \tau^2}\bigg)^{Y_i+\alpha}
   t^{-\alpha-a-1} L(t) dt \\    
    &\leq (\tau^2)^{-Y_i-\alpha}
		\int_{\frac{1}{\tau^2}(\frac{1}{\epsilon}-1)} ^{\infty} t^{-\alpha-a-1} L(t) dt \hspace{0.05cm}. \tag{25}
\end{align*}
	Next using Lemma \ref{chap4lem1} and the definition of $L(\cdot)$, we obtain 
	\begin{equation*} \label{eq:T-3.3}
		\int_{\frac{1}{\tau^2}(\frac{1}{\epsilon}-1)} ^{\infty} t^{-\alpha-a-1} L(t) dt= \frac{1}{(a+\alpha)} (\tau^2)^{a+\alpha} \bigg(\frac{\epsilon}{1-\epsilon} \bigg)^{a+\alpha} L \bigg(\frac{1}{\tau^2}\bigg)  (1+o(1)) \hspace{0.05cm}. \tag{26}
	\end{equation*}
	Using \eqref{eq:T-3.1}-\eqref{eq:T-3.3} the proof of Theorem \ref{sparsecountshrinktun1} is obtained easily. 
\end{proof}

\begin{proof}[Proof of Theorem \ref{sparsecountshrinktun2}]
	Using the definition,
	\begin{equation*}
		\mathbb{E}(1-\kappa_i|Y_i, \tau) = \frac{\int_{0}^{1}	{\kappa_i^{a+\alpha-1}}(1-\kappa_i)^{Y_i-a}L\bigg(\frac{1}{\tau^2}(\frac{1}{\kappa_i}-1)\bigg)d\kappa_i}{\int_{0}^{1}	\kappa_i^{a+\alpha-1}(1-\kappa_i)^{Y_i-a-1}L\bigg(\frac{1}{\tau^2}(\frac{1}{\kappa_i}-1)\bigg)d\kappa_i} \hspace{0.05cm}.
	\end{equation*}
	Again using the change of variable $t=\frac{1}{\tau^2}(\frac{1}{\kappa_i}-1)$ in the numerator and denominator of the previous quantity and applying \eqref{chap4eq:L-2.2} of Lemma \ref{chap4lem2} for the denominator, we have
	\begin{equation*} \label{eq:T-4.1}
		\mathbb{E}(1-\kappa_i|Y_i, \tau) \leq \frac{a}{c_0} (K_0^{-a}-K_1^{-a})^{-1} \tau^2 (1+K_1 \tau^2)^{(Y_i+\alpha)} \int_{0}^{\infty} (1+t \tau^2)^{-(Y_i+\alpha+1)} t^{Y_i-a} L(t) dt \hspace{0.05cm}. \tag{27}
	\end{equation*}
    Next, we divide the range of integration into two parts, namely in $t \in (0,1)$ and $t \geq 1$. We observe that
	\begin{equation*} \label{eq:T-4.2}
		\tau^2 \int_{0}^{1} (1+t \tau^2)^{-(Y_i+\alpha+1)} t^{Y_i-a} L(t) dt  \leq \tau^2 K^{-1} \hspace{0.05cm}. \tag{28}
	\end{equation*}
	Next using \hyperlink{assumption2}{Assumption 2} on $L(\cdot)$, it is easy to see that for any $Y_i \in [0, a-1)$, 
	\begin{equation*} \label{eq:T-4.3}
		\tau^2 \int_{1}^{\infty} (1+t \tau^2)^{-(Y_i+\alpha+1)} t^{Y_i-a} L(t) dt  \leq \frac{\tau^2 M}{(a-Y_i-1)} \hspace{0.05cm}. \tag{29}
	\end{equation*}
	Finally, Theorem \ref{sparsecountshrinktun2} follows from \eqref{eq:T-4.1}-\eqref{eq:T-4.3}.
\end{proof}

\begin{proof}[Proof of Theorem \ref{sparsecountshrinktun3}]
	Fix any $\eta \in (0,1)$ and $\delta_1 \in (0,1)$. Now using the definition,
	\begin{equation*}
		\mathbb{P}(\kappa_i> \eta|Y_i, \tau)= \frac{\int_{\eta}^{1} \kappa_i^{a+\alpha-1}(1-\kappa_i)^{Y_i-a-1}L\bigg(\frac{1}{\tau^2}(\frac{1}{\kappa_i}-1)\bigg)d\kappa_i}{\int_{0}^{1}	\kappa_i^{a+\alpha-1}(1-\kappa_i)^{Y_i-a-1}L\bigg(\frac{1}{\tau^2}(\frac{1}{\kappa_i}-1)\bigg)d\kappa_i} \hspace{0.05cm}.
	\end{equation*}
	Then using the change of variable  $t= \frac{1}{\tau^2}(\frac{1}{\kappa_i}-1)$ to  both the numerator and denominator of the right-hand side of the above equation, we obtain, for any $\delta_1 \in (0,1)$,
	\begin{equation*}
		\mathbb{P}(\kappa_i> \eta|Y_i, \tau)= \frac{\int_{0}^{\frac{1}{\tau^2}(\frac{1}{\eta}-1)}(1+t \tau^2)^{-(Y_i+\alpha)} t^{Y_i-a-1} L(t) dt}{\int_{0}^{\infty}(1+t \tau^2)^{-(Y_i+\alpha)} t^{Y_i-a-1} L(t) dt}
	\end{equation*}
	\begin{equation*} \label{eq:T-5.1}
		\leq \frac{\int_{0}^{\frac{1}{\tau^2}(\frac{1}{\eta}-1)}(1+t \tau^2)^{-(Y_i+\alpha)} t^{Y_i-a-1} L(t) dt}{\int_{\frac{1}{\tau^2}(\frac{1}{\eta \delta_1}-1)}^{\infty}(1+t \tau^2)^{-(Y_i+\alpha)} t^{Y_i-a-1} L(t) dt} \hspace{0.05cm}. \tag{30}
	\end{equation*}
	Note that the numerator of \eqref{eq:T-5.1} can be further bounded as
	\begin{align*} \label{eq:T-5.2}
		\int_{0}^{\frac{1}{\tau^2}(\frac{1}{\eta}-1)}(1+t \tau^2)^{-(Y_i+\alpha)} t^{Y_i-a-1} L(t) dt & = (\tau^2)^{-Y_i} \int_{0}^{\frac{1}{\tau^2}(\frac{1}{\eta}-1)} (\frac{t \tau^2}{1+t \tau^2})^{Y_i} (1+t \tau^2)^{-\alpha} t^{-a-1} L(t) dt \\
		& \leq (\tau^2)^{-Y_i} (1-\eta)^{Y_i} \int_{0}^{\infty}t^{-a-1}L(t) dt\\
		&= K^{-1} \bigg(\frac{1-\eta}{\tau^2} \bigg)^{Y_i} \hspace{0.05cm}. \tag{31}
	\end{align*}
	Here the inequality occurs due to the fact that $t\tau^2/(1+t \tau^2)$ is increasing in $t$ for any fixed $\tau>0$ whenever $t \in (0,\frac{1}{\tau^2}(\frac{1}{\eta}-1))$. Next we use $\int_{0}^{\infty}t^{-a-1}L(t) dt=K^{-1}$. \\
	On the other hand, for the denominator, we have
	\begin{align*} \label{eq:T-5.3}
		& \int_{\frac{1}{\tau^2}(\frac{1}{\eta \delta_1}-1)}^{\infty}(1+t \tau^2)^{-(Y_i+\alpha)} t^{Y_i-a-1} L(t) dt \\
		& = (\tau^2)^{-Y_i-\alpha} \int_{\frac{1}{\tau^2}(\frac{1}{\eta \delta_1}-1)}^{\infty} \bigg(\frac{t \tau^2}{1+t \tau^2}\bigg)^{Y_i+\alpha} t^{-a-\alpha-1}L(t) dt \\
		& \geq (\tau^2)^{-Y_i-\alpha} (1-\eta \delta_1)^{Y_i+\alpha}  \int_{\frac{1}{\tau^2}(\frac{1}{\eta \delta_1}-1)}^{\infty} t^{-a-\alpha-1}L(t) dt \\
		& \geq c_0 (\tau^2)^{-Y_i-\alpha} (1-\eta \delta_1)^{Y_i+\alpha} \int_{\frac{1}{\tau^2}(\frac{1}{\eta \delta_1}-1)}^{\infty} t^{-a-\alpha-1} dt \\
		&= \frac{c_0}{(a+\alpha)} \bigg(\frac{1}{\tau^2}\bigg)^{Y_i-a} (1-\eta \delta_1)^{Y_i-a} (\eta \delta_1)^{a+\alpha}  \hspace{0.05cm}. \tag{32}
	\end{align*}
	In the chain of inequalities, the first one holds due to the fact that $t\tau^2/(1+t \tau^2)$ is increasing in $t$ for any fixed $\tau>0$ whenever $t \geq \frac{1}{\tau^2}(\frac{1}{\eta \delta_1}-1)$. The second inequality follows from \hyperlink{assumption2}{Assumption 2} on $L(\cdot)$ by noting that for any fixed $\eta \in (0,1)$ and $\delta_1 \in (0,1)$, $\tau$ can be made small enough so that $\frac{1}{\tau^2}(\frac{1}{\eta \delta_1}-1) \geq t_0$. Finally, \eqref{eq:T-5.2} and \eqref{eq:T-5.3} provide the upper bound on \eqref{eq:T-5.1} and complete the proof of Theorem \ref{sparsecountshrinktun3}.
\end{proof}

\begin{proof}[Proof of Theorem \ref{sparsecounttypeItun}]
    In order to calculate type I error, note that, the posterior distribution of $\kappa_i$ given $Y_i$ and $\tau$ has the same from for all $i=1,2,\cdots,n$. Also, under $H_{0i}$, the marginal distribution of $Y_i$ is the same for all $i=1,2,\cdots,n$. These two facts indicate that the probability of type I error of $i^{\text{th}}$ decision rule, denoted as $t_{1i}$ is independent of $i$.
    We further then have,
	\begin{align*} \label{eq:T-7.2}
	t_{1i} \equiv	t_1 &= \mathbb{P}_{H_{0i}} (\mathbb{E}(1-\kappa_i|Y_i, \tau)> \frac{\delta}{2(\delta+1)}) \\
		&\leq  \mathbb{P}_{H_{0i}} (\mathbb{E}(1-\kappa_i|Y_i, \tau)> \frac{\delta}{2(\delta+1)}, Y_i \leq \frac{a}{2}) + \mathbb{P}_{H_{0i}} (Y_i > \frac{a}{2}) \hspace{0.05cm}. \tag{33}
	\end{align*}
	Now concentrate on the first term in \eqref{eq:T-7.2}. In this case, 
    using Theorem \ref{sparsecountshrinktun2}, as $\tau \to 0$, this probability can be further bounded by
	\begin{align*} \label{eq:T-7.3}
       & \mathbb{P}_{H_{0i}} \bigg(M_1 \tau^2[K^{-1}+\frac{M}{\frac{a}{2}-1}]  (1+K_1\tau^2)^{(Y_i+\alpha)}  > \frac{\delta}{2(\delta+1)} , Y_i \leq \frac{a}{2} \bigg ) \\ 
        & \leq 
 \mathbb{P} \bigg(M_2 \tau^2 (1+K_1\tau^2)^{(\frac{a}{2}+\alpha)}  > \frac{\delta}{2(\delta+1)}  \bigg ) \\
    &= 0, \textrm{for all sufficiently large } n,    \tag{34} 
	\end{align*}
    where $M_1$ and $M_2$ are constants depending only on $c_0, a$, $\alpha, K_0$ and $K_1$ and are independent of $Y_i$.
	Since, under $H_{0i}$, $Y_i \simiid NB(\alpha, \frac{1}{\beta+1})$, for the second term of \eqref{eq:T-7.2}, using the Markov's inequality, we have
	\begin{equation*} \label{chap4eq:T-7.4}
		\mathbb{P}_{H_{0i}} (Y_i > \frac{a}{2}) \leq \frac{2\alpha \beta}{a} \hspace{0.05cm}. \tag{35}
	\end{equation*}
	The proof is completed using \eqref{eq:T-7.2}-\eqref{chap4eq:T-7.4}. 
\end{proof}
\begin{proof}[Proof of Theorem \ref{sparsecounttypeIItun}]
    Note that
    under $H_{1i}$, the distribution of $Y_i$ is the same for all $i$, $i=1,2,\cdots,n$. Also, note that, the form of the posterior distribution of $\kappa_i$ given $Y_i$ and $\tau$ is same for all $i=1,2,\cdots,n$. As a consequence of this, the probability of type II error of $i^{\text{th}}$ decision rule is denoted as  $t_2$. \\
Hence, as $n \to \infty$, the upper bound on $t_2$ is of the form
	\begin{align*} \label{chap4eq:T-8.2}
	\lim_{n \to \infty}	t_2 &= \lim_{n \to \infty} \mathbb{P}_{H_{1i}}(\mathbb{E}(\kappa_i|Y_i,\tau)\geq\frac{\delta+2}{2(\delta+1)}) \\
		& \leq \lim_{n \to \infty}\mathbb{P}_{H_{1i}}(\mathbb{E}(\kappa_i|Y_i,\tau)\geq\frac{\delta+2}{2(\delta+1)}, Y_i >a) + \lim_{n \to \infty}\mathbb{P}_{H_{1i}} (Y_i \leq a).
        \tag{36}
	\end{align*}
First note that
\begin{align*} \label{chap4typeIInew1}
  \lim_{n \to \infty}\mathbb{P}_{H_{1i}} (Y_i \leq a)   & =  \lim_{n \to \infty} \sum_{y=0}^{[a]} \binom{y+\alpha-1}{y} (1-\frac{1}{\beta+\delta+1})^y (\beta+\delta+1)^{-\alpha}\\
  &= \sum_{y=0}^{[a]} \binom{y+\alpha-1}{y} (1-\frac{1}{\delta+1})^y (\delta+1)^{-\alpha}= \mathbb{P}(Y \leq a) \tag{37},
\end{align*}
where $Y \sim NB(\alpha, \frac{1}{\delta+1})$ and $[x]$ denotes the greatest integer less than or equal to $x$.
Note that
\begin{align*}\label{chap4eq:T-8.3}
    \lim_{n \to \infty}\mathbb{P}_{H_{1i}}(\mathbb{E}(\kappa_i|Y_i,\tau)\geq\frac{\delta+2}{2(\delta+1)}, Y_i >a) = \lim_{n \to \infty} \mathbb{E}_{H_{1i}}(Z_n),
\end{align*}
where $Z_n = 1\{\mathbb{E}(\kappa_i|Y_i,\tau)\geq\frac{\delta+2}{2(\delta+1)}, Y_i >a\}$. 
Therefore
\begin{align*}
    \mathbb{E}_{H_{1i}}(Z_n) &= \int_{\Omega}1\{\mathbb{E}(\kappa|y,\tau)\geq\frac{\delta+2}{2(\delta+1)}, y >a\} f_n(y) d\mu\\
    &= \int_{\Omega} g_n(y) d\mu, \textrm{say},
\end{align*}
where $f_n(y)$ denotes the probability mass function of $Y_i$ at $y$ under $H_{1i}$ and is of the form
\begin{align*}
    f_n(y) & = \binom{y+\alpha-1}{y} (1-\frac{1}{\beta+\delta+1})^y (\beta+\delta+1)^{-\alpha}, y=0,1,2,\cdots
\end{align*}
and $\mu$ is the counting measure on $\Omega =\{0,1,2,\cdots\}$.
Here,
\begin{equation*} \label{chap4eq:T-8.4}
    {E}(\kappa|y, \tau) = \frac{\int_{0}^{1} \kappa_i^{a+\alpha}(1-\kappa)^{y-a-1}L\bigg(\frac{1}{\tau^2}(\frac{1}{\kappa}-1)\bigg)d\kappa}{\int_{0}^{1}	\kappa^{a+\alpha-1}(1-\kappa)^{y-a-1}L\bigg(\frac{1}{\tau^2}(\frac{1}{\kappa}-1)\bigg)d\kappa} \hspace{0.05cm}. \tag{38}
\end{equation*}
We now show that $f_n(y) \leq h(y)$ for all $y \in \Omega$, for some $h(y)$ such that it's distribution is independent of $n$ and $\int_{\Omega}h(y) d \mu < \infty$. Since, $\beta \to 0$ as $n \to \infty$, for sufficiently large $n$, we assume $\beta \leq 1$. As a result,
\begin{align*}
    f_n(y) & \leq \binom{y+\alpha-1}{y} (1-\frac{1}{\delta+2})^y (\delta+1)^{-\alpha} \\
    &= \binom{y+\alpha-1}{y} (1-\frac{1}{\delta+1})^y (\delta+1)^{-\alpha} \frac{(1-\frac{1}{\delta+2})^y}{(1-\frac{1}{\delta+1})^y} \\
    &= \binom{y+\alpha-1}{y} (1-\frac{1}{\delta+1})^y (\delta+1)^{-\alpha}  (1+\frac{1}{\delta^2+2\delta})^y.
\end{align*}
 Hence, if we define $h(y)= \binom{y+\alpha-1}{y} (1-\frac{1}{\delta+1})^y (\delta+1)^{-\alpha}  (1+\frac{1}{\delta^2+2\delta})^y$, then 
\begin{align*}
    \sum_{y=0}^{\infty} h(y) &= \mathbb{E}[(1+\frac{1}{\delta^2+2\delta})^Y]= \mathbb{E}[e^{Y\log (1+\frac{1}{\delta^2+2\delta})}],
\end{align*}
where $Y \sim NB(\alpha,\frac{1}{\delta+1})$. Recall that if $X \sim NB(r,q)$, then the moment generating function (MGF) of $t$ exists if $t<-\log(1-q)$. Here $q=\frac{1}{\delta+1}$. Therefore
$\mathbb{E}[e^{Y\log (1+\frac{1}{\delta^2+2\delta})}]$ exists if $\log (1+\frac{1}{\delta^2+2\delta}) < \log (1+\frac{1}{\delta})$, i.e.,
for any $\delta>0$. Which implies, $f_n(y) \leq h(y)= \binom{y+\alpha-1}{y} (1-\frac{1}{\delta+1})^y (\delta+1)^{-\alpha}  (1+\frac{1}{\delta^2+2\delta})^y$ and $\int_{\Omega}h(y) d \mu < \infty$.
Also, note that, for all $y \in \Omega, \lim_{n \to \infty} f_n(y)= f(y)=\binom{y+\alpha-1}{y} (1-\frac{1}{\delta+1})^y (\delta+1)^{-\alpha}$. Define, $m(y)= \lim_{\tau \to 0} \mathbb{E}(\kappa|y,\tau)$.
Hence,
 for all $y \in \Omega$,
\begin{equation*}
     \lim_{n \to \infty} g_n(y) = 1\{m(y)\geq\frac{\delta+2}{2(\delta+1)}, y >a\}\binom{y+\alpha-1}{y} (1-\frac{1}{\delta+1})^y (\delta+1)^{-\alpha}= g(y), \textrm{say}.
\end{equation*}
Therefore, by the Dominated Convergence Theorem,
\begin{equation*}
    \int_{\Omega} g_n(y) d \mu \to \int_{\Omega} g(y) d\mu \textrm{ as } n \to \infty.
\end{equation*}

    
Under \hyperlink{assumption2}{Assumption 2} of Section \ref{sec-4.3}, for $y >a$ using \eqref{chap4eq:T-8.4}, we have
\begin{equation*}
    m(y) = \frac{Beta(a+\alpha+1, y-a)}{Beta(a+\alpha, y-a)}=\frac{a+\alpha}{y+\alpha}.
\end{equation*}
Hence, we obtain
\begin{align*} \label{chap4eq:T-8.5}
    \lim_{n \to \infty}\mathbb{P}_{H_{1i}}(\mathbb{E}(\kappa_i|Y_i,\tau)\geq\frac{\delta+2}{2(\delta+1)}, Y_i >a) 
   & =  \mathbb{P}(m(Y)\geq\frac{\delta+2}{2(\delta+1)}, Y >a) \\
    &= \mathbb{P} ( \frac{a+\alpha}{Y+\alpha}\geq\frac{\delta+2}{2(\delta+1)}, Y >a) \\
    &=
    \mathbb{P} \bigg(a<Y \leq 2a+\alpha -\frac{2(a+\alpha)}{(\delta+2)}\bigg),
    \tag{39}
\end{align*}
where $Y \sim NB(\alpha,\frac{1}{\delta+1})$.
Finally the proof of Theorem \ref{sparsecounttypeIItun} is obtained by combining \eqref{chap4eq:T-8.2}, \eqref{chap4typeIInew1} and \eqref{chap4eq:T-8.5}.

\end{proof}

\begin{proof}[Proof of Theorem \ref{sparsecounttypeIEB}]
	First, let us find an expression for $\gamma_n= P(Y_i \geq 1)$. We have,
	\begin{align*}
		\gamma_n &= (1-p) \mathbb{P}_{H_{0i}}(Y_i \geq 1)+p  \mathbb{P}_{H_{1i}}(Y_i \geq 1) \\
		&= p[\mathbb{P}_{H_{1i}}(Y_i \geq 1)+ \frac{(1-p)}{p} \mathbb{P}_{H_{0i}}(Y_i \geq 1)]
	\end{align*}
	Since, $Y_i \simiid NB(\alpha, \frac{1}{\beta+1})$,  under $H_{0i}$, using Markov's inequality and \hyperlink{assumption1}{Assumption 1}, we have, $\frac{1}{p} \mathbb{P}_{H_{0i}}(Y_i \geq 1) \to 0$ as $n \to \infty$. Also, note that, under $H_{1i}$,
	\begin{equation*}
		\mathbb{P}(Y_i \geq 1)= (1-(\beta+\delta+1)^{-\alpha}).
	\end{equation*}
	Combining these two, we obtain, as $n \to \infty$,
    \begin{equation*} \label{eq:newT-9.1}
        \gamma_n= (1-(\beta+\delta+1)^{-\alpha})p(1+o(1)). \tag{40}
    \end{equation*}
	Now, using the definition of type I error of the $i^{\text{th}}$ decision rule in \eqref{eq:4.3.8}, we have
	\begin{align*}\label{eq:T-9.1}
		t^{\text{EB}}_{1i} &= \mathbb{P}_{H_{0i}} (\mathbb{E}(1-\kappa_i|Y_i, \widehat{\tau})> \frac{\delta}{2(\delta+1)}) \\
		&= \mathbb{P}_{H_{0i}} (\mathbb{E}(1-\kappa_i|Y_i, \widehat{\tau})> \frac{\delta}{2(\delta+1)}, 0 \leq Y_i \leq \frac{a}{2})+\mathbb{P}_{H_{0i}} (\mathbb{E}(1-\kappa_i|Y_i, \widehat{\tau})> \frac{\delta}{2(\delta+1)}, Y_i > \frac{a}{2}) \hspace{0.05cm}. \tag{41}
	\end{align*}
	For the second term of \eqref{eq:T-9.1}, using Markov's inequality, we have
	\begin{equation*} \label{eq:T-9.2}
		\mathbb{P}_{H_{0i}} (\mathbb{E}(1-\kappa_i|Y_i, \widehat{\tau})> \frac{\delta}{2(\delta+1)}, Y_i > \frac{a}{2}) \leq   \mathbb{P}_{H_{0i}} (Y_i > \frac{a}{2}) \leq \frac{2\alpha \beta}{a} \hspace{0.05cm}. \tag{42}
	\end{equation*}
	For the first term in \eqref{eq:T-9.1}, we obtain
	\begin{align*} \label{eq:T-9.3}
		&  \mathbb{P}_{H_{0i}} (\mathbb{E}(1-\kappa_i|Y_i, \widehat{\tau})> \frac{\delta}{2(\delta+1)}, 0 \leq Y_i \leq \frac{a}{2}) \\
		&=  \mathbb{P}_{H_{0i}} (\mathbb{E}(1-\kappa_i|Y_i, \widehat{\tau})> \frac{\delta}{2(\delta+1)}, 0 \leq Y_i \leq \frac{a}{2}, \widehat{\tau} \leq 2 \gamma_n) \\ &+  \mathbb{P}_{H_{0i}} (\mathbb{E}(1-\kappa_i|Y_i, \widehat{\tau})> \frac{\delta}{2(\delta+1)}, 0 \leq Y_i \leq \frac{a}{2}, \widehat{\tau} > 2 \gamma_n) \hspace{0.05cm}. \tag{43}
	\end{align*}
	Now, using the form of $\gamma_n$,
   we note that 
	\begin{align*} \label{eq:T-9.4}
		&  \mathbb{P}_{H_{0i}} (\mathbb{E}(1-\kappa_i|Y_i, \widehat{\tau})> \frac{\delta}{2(\delta+1)}, 0 \leq Y_i \leq \frac{a}{2}, \widehat{\tau} \leq 2 \gamma_n) \\
		& \leq \mathbb{P}_{H_{0i}} (\mathbb{E}(1-\kappa_i|Y_i,  2 \gamma_n)> \frac{\delta}{2(\delta+1)}, 0 \leq Y_i \leq \frac{a}{2}) \\
        &= 0, \textrm{ for all sufficiently large } n.
        \tag{44}
	\end{align*}
    The inequality above holds due to the fact that for any fixed $y \geq 0$, $\mathbb{E}(1-\kappa|y,\tau)$ is non-decreasing in $\tau$. The last assertion follows from \eqref{eq:T-7.3}. \\
	\hspace*{0.05cm} Now let us concentrate on the second term of \eqref{eq:T-9.3}. Define $\widehat{\tau}_{1}= \frac{1}{n}$ and $\widehat{\tau}_{2}= \frac{1}{n} \sum_{i=1}^{n} \mathbf{1} \{Y_i \geq 1 \}$. Hence, $\widehat{\tau}= \text{max} \{\widehat{\tau}_1,\widehat{\tau}_2\}$. Also, note that $\gamma_n \sim (1-(\beta+\delta+1)^{-\alpha})p$ and $p \propto n^{-\epsilon}, 0<\epsilon<1$ implies $\frac{1}{n} < 2 \gamma_n$ for all sufficiently large $n$. Next, we use the fact that $\{\widehat{\tau} > 2 \gamma_n\}$ implies that either of the two events $\{\widehat{\tau}_1 > 2 \gamma_n\}$ or $\{\widehat{\tau}_2 > 2 \gamma_n\}$ will happen. Using these facts for all sufficiently large $n$, we have
	\begin{align*} \label{eq:T-9.5}
		& \mathbb{P}_{H_{0i}} (\mathbb{E}(1-\kappa_i|Y_i, \widehat{\tau})> \frac{\delta}{2(\delta+1)}, 0 \leq Y_i \leq \frac{a}{2}, \widehat{\tau} > 2 \gamma_n)  \leq \mathbb{P}_{H_{0i}}(\widehat{\tau} > 2 \gamma_n) \\
		& \leq  \mathbb{P}_{H_{0i}}(\widehat{\tau}_1 > 2 \gamma_n)+  \mathbb{P}_{H_{0i}}(\widehat{\tau}_2 > 2 \gamma_n) \\
		&= \mathbb{P}_{H_{0i}}(\widehat{\tau}_2 > 2 \gamma_n) \\
		& \leq \mathbb{P}_{H_{0i}}(\widehat{\tau}_2 > 2 \gamma_n, Y_i < 1)+ \mathbb{P}_{H_{0i}}(Y_i \geq 1) \\
		& \leq \mathbb{P}_{H_{0i}}(\widehat{\tau}_2 > 2 \gamma_n, Y_i < 1) +{\alpha \beta}{} \hspace{0.05cm}.  \tag{45}
	\end{align*}
	Now we are only left with only the first term of \eqref{eq:T-9.5}. Note that
	\begin{equation*}
		\{\widehat{\tau}_2 > 2 \gamma_n, Y_i < 1\} \subseteq \{ \frac{1}{n} \sum_{\substack{j=1 \\ (j \neq i)}}^{n} \mathbf{1} \{Y_j \geq 1 \} >2 \gamma_n \}.
	\end{equation*}
	As a result, using the Chernoff-Hoeffding inequality, we have
	\begin{align*} \label{eq:T-9.6}
		\mathbb{P}_{H_{0i}}(\widehat{\tau}_2 > 2 \gamma_n, Y_i < 1) & \leq \mathbb{P}_{H_{0i}}(\frac{1}{n} \sum_{\substack{j=1 \\ (j \neq i)}}^{n} \mathbf{1} \{Y_j \geq 1\} >2 \gamma_n )\\
		& \leq \mathbb{P}(\frac{1}{n-1} \sum_{\substack{j=1 \\ (j \neq i)}}^{n} \mathbf{1} \{Y_j \geq 1\} \geq 2 \gamma_n ) \\
		& \leq e^{-(n-1)D(2  \gamma_n || \gamma_n)} \hspace{0.05cm},  \tag{46}
	\end{align*}
	where $D(x || y)$ is defined as, $D(x || y)= x \log (\frac{x}{y})+(1-x) \log(\frac{1-x}{1-y})$. Next using these facts that $\log (\frac{1}{1-x}) \sim x$ as $x \to 0$ and $\gamma_n \to 0$ as $n \to \infty$ and also using the calculations of Theorem 10 of \cite{ghosh2016asymptotic},
    we have
	\begin{equation*} \label{eq:T-9.7}
		D(2  \gamma_n || \gamma_n)= (2 \log 2-1)\gamma_n (1+o(1))=(2 \log 2-1)(1-(\beta+\delta+1)^{-\alpha})p(1+o(1)) \hspace{0.05cm}.  \tag{47}
	\end{equation*}
	Here $o(1)$ is non-random and 
    tends to $0$ as $n \to \infty$. Hence using \eqref{eq:T-9.6} and \eqref{eq:T-9.7}, we obtain
	\begin{equation*} \label{eq:T-9.8}
		\mathbb{P}_{H_{0i}}(\widehat{\tau}_2 > 2 \gamma_n, Y_i < 1) \leq e^{-(2 \log 2-1)(1-(\beta+\delta+1)^{-\alpha})np(1+o(1))} \hspace{0.05cm}.  \tag{48}
	\end{equation*}
	Combining the above arguments the proof of Theorem \ref{sparsecounttypeIEB} follows.
\end{proof}

\begin{proof}[Proof of Theorem \ref{sparsecounttypeIIEB}]
	Let us fix any $C_3 \in (0,1)$. Now, using the definition of type II error of the $i^{\text{th}}$ decision rule in \eqref{eq:4.3.8}, we have
	\begin{align*} \label{eq:T-10.1}
		t^{\text{EB}}_{2i} &= \mathbb{P}_{H_{1i}} (\mathbb{E}(\kappa_i|Y_i, \widehat{\tau}) \geq \frac{\delta+2}{2(\delta+1)}) \\
		&= \mathbb{P}_{H_{1i}} (\mathbb{E}(\kappa_i|Y_i, \widehat{\tau})\geq \frac{\delta+2}{2(\delta+1)}, \widehat{\tau} < C_3 \gamma_n)+  \mathbb{P}_{H_{1i}} (\mathbb{E}(\kappa_i|Y_i, \widehat{\tau})\geq \frac{\delta+2}{2(\delta+1)}, \widehat{\tau} \geq C_3 \gamma_n)
		\hspace{0.05cm}, \tag{49}
	\end{align*}
    where $\gamma_n$ is same as defined in the previous theorem. Let us first bound the second term in the r.h.s. of \eqref{eq:T-10.1}.
	Using the fact that for any fixed $y$, $\mathbb{E}(\kappa|y,\tau)$ is decreasing in $\tau$, towards that, we first note that
	\begin{equation*}  \label{eq:T-10.2}
		\{\mathbb{E}(\kappa_i|Y_i, \widehat{\tau})\geq \frac{\delta+2}{2(\delta+1)}, \widehat{\tau} \geq C_3 \gamma_n\} \subseteq \{\mathbb{E}(\kappa_i|Y_i,C_3 \gamma_n) \geq\frac{\delta+2}{2(\delta+1)}\}  \hspace{0.05cm}. \tag{50}
	\end{equation*}
	Now applying the same set of arguments used in the proof of Theorem \ref{sparsecounttypeIItun} and noting under ${H_{1i}}$, $Y_i \simiid NB(\alpha, \frac{1}{\beta+ \delta+1})$, we obtain, as $n \to \infty$,
	\begin{align*}  \label{eq:T-10.3}
		\mathbb{P}_{H_{1i}}(\mathbb{E}(\kappa_i|Y_i,C_3 \gamma_n)\geq \frac{\delta+2}{2(\delta+1)}) & \leq 
      \mathbb{P} \bigg(Y \leq 2a+\alpha -\frac{2(a+\alpha)}{(\delta+2)}\bigg),
        \tag{51}
	\end{align*}
    where $Y \sim NB(\alpha,\frac{1}{\delta+1})$.
    Now our aim is to show that the first term in \eqref{eq:T-10.1} goes to $0$ as $n \to \infty$. Note that, $\widehat{\tau}_2 \geq \frac{1}{n} \sum_{\substack{j=1 \\ (j \neq i)}}^{n} \mathbf{1} \{Y_j \geq 1 \}$. Using this observation and noting that due to independence, the distribution of remaining $Y_j$'s do not depend on that of $Y_i$, we have
	\begin{align*}\label{eq:T-10.4}
		\mathbb{P}_{H_{1i}} (\mathbb{E}(\kappa_i|Y_i, \widehat{\tau})\geq
		\frac{1}{2}, \widehat{\tau} < C_3 \gamma_n) & \leq \mathbb{P}_{H_{1i}}(\widehat{\tau} < C_3 \gamma_n) \\
		& \leq \mathbb{P}_{H_{1i}}(\widehat{\tau}_2 < C_3 \gamma_n) \\
		& \leq \mathbb{P}_{H_{1i}}(\frac{1}{n} \sum_{\substack{j=1 \\ (j \neq i)}}^{n} \mathbf{1} \{Y_j \geq 1 \} \leq C_3 \gamma_n) \\
		&= \mathbb{P}(-\frac{1}{n-1}(\sum_{\substack{j=1 \\ (j \neq i)}}^{n}( \mathbf{1} \{Y_j \geq 1 \}-\gamma_n)) \geq (1-\frac{n}{n-1}C_3)\gamma_n)\hspace{0.05cm}, \tag{52}
	\end{align*}
    where probability in the last statement is calculated using the marginal distribution \eqref{eq:4.1.2} of $Y_j, j \neq i$.
	Since, $1-\frac{n}{n-1}C_3 \to 1-C_3$ as $n \to \infty$ and $C_3 \in (0,1)$, so, $1-\frac{n}{n-1}C_3 >0$ for all sufficiently large $n$. Finally using Markov's inequality,
	\begin{align*} \label{eq:T-10.5}
		&  \mathbb{P}(-\frac{1}{n-1}(\sum_{\substack{j=1 \\ (j \neq i)}}^{n} (\mathbf{1} \{Y_i \geq 1 \}-\gamma_n)) \geq (1-\frac{n}{n-1}C_3)\gamma_n) \\
		&\leq \mathbb{P}\bigg(|\frac{1}{n-1}(\sum_{\substack{j=1 \\ (j \neq i)}}^{n} (\mathbf{1} \{Y_i \geq 1 \}-\gamma_n))| \geq (1-\frac{n}{n-1}C_3)\gamma_n \bigg) \\
		& \leq \frac{\gamma_n(1-\gamma_n)}{(n-1)(1-\frac{n}{n-1}C_3)^2\gamma^2_n} \\
		&= \frac{(1-\gamma_n)}{(1-C_3)^2 n \gamma_n}(1+o(1))=o(1) \hspace{0.05cm}. \tag{53}
	\end{align*}
	At the final step, we use $\gamma_n \sim (1-(\beta+\delta+1)^{-\alpha})p$ and $p \propto n^{-\epsilon}, 0< \epsilon<1$. Combining \eqref{eq:T-10.4} and \eqref{eq:T-10.5}, we obtain, as $n \to \infty$,
	\begin{equation*}
		\mathbb{P}_{H_{1i}} (\mathbb{E}(\kappa_i|Y_i, \widehat{\tau})> \frac{\delta+2}{2(\delta+1)}, \widehat{\tau} < C_3 \gamma_n)=o(1) \hspace{0.05cm}.
	\end{equation*}
	Hence using \eqref{eq:T-10.3}, as $n \to \infty$, we have
	\begin{equation*} \label{eq:T-10.6}
		t^{\text{EB}}_{2i} \leq
        \mathbb{P} \bigg(Y \leq 2a+\alpha -\frac{2(a+\alpha)}{(\delta+2)}\bigg)
       (1+o(1)) \hspace{0.05cm}, \tag{54}  
	\end{equation*}
    since the r.h.s. of \eqref{eq:T-10.3} is bounded away from zero under \hyperlink{assumption1}{Assumption 1}.
    This completes the proof of Theorem \ref{sparsecounttypeIIEB}.
\end{proof}

\begin{proof}[Proof of Theorem \ref{sparsecountabostun}]
	First, we study the probabilities of type I error ($t^{\text{BO}}_1$) and type II error ($t^{\text{BO}}_2$) of the decision rule when each $\theta_i$ is modeled by a two-group prior of the form \eqref{eq:4.1.1}. \\
	Since, under $H_{0i}$, $Y_i \simiid NB(\alpha, \frac{1}{\beta+1})$, using Markov's inequality and \hyperlink{assumption1}{Assumption 1}, $t^{\text{BO}}_1$ can be bounded as
	\begin{align*}  \label{chap4eq:T-1.1}
		t^{\text{BO}}_1 &= \mathbb{P}_{H_{0i}}(Y_i > C_{p,\alpha,\beta, \delta})\\
		&= \mathbb{P}_{H_{0i}} \bigg(Y_i > \frac{\log(\frac{1}{p} ) +\alpha \log (1+\delta)}{\log (\frac{1}{\beta})}(1+o(1))\bigg) \\
        & = \mathbb{P}_{H_{0i}} (Y_i > \frac{1}{C_1}(1+o(1))) \\
		& \leq \alpha \beta C_1 (1+o(1))  \hspace{0.05cm}. \tag{55}
	\end{align*}
    Here, inequality in \eqref{chap4eq:T-1.1} holds due to Markov's inequality.
	Since under the assumption $\beta \propto p^{C_1}$, for some $C_1>1$, we have, $\frac{t^{\text{BO}}_1}{p} \to 0$ as $n \to \infty$. Again using \hyperlink{assumption1}{Assumption 1}, we have
	\begin{align*}  \label{chap4eq:T-1.2}
		t^{\text{BO}}_2 &= \mathbb{P}_{H_{1i}}(Y_i \leq  C_{p,\alpha,\beta, \delta})\\
		&=  \mathbb{P}_{H_{1i}} \bigg(Y_i \leq \frac{1}{C_1}(1+o(1))\bigg).
        \tag{56}
	\end{align*}
    Since under $H_{1i}, Y_i \simiid NB (\alpha, \frac{1}{\beta+\delta+1})$, using the probability mass function (p.m.f.) of a negative binomial random variable, we have
    \begin{align*}
       t^{\text{BO}}_2 & =  \mathbb{P}_{H_{1i}} (Y_i =0) 
        = (\beta+\delta+1)^{-\alpha}.
    \end{align*}
    Hence under \hyperlink{assumption1}{Assumption 1}, 
 we have, for all sufficiently large $n$,
 \begin{align*} \label{chapbotypeIIlower}
     t^{\text{BO}}_2 & = (\delta+1)^{-\alpha}(1+o(1)). \tag{57}
 \end{align*}
    As a consequence of this, the asymptotic expression
	for the optimal Bayes risk based on the decision rule \eqref{eq:4.2.3} with \eqref{eq:4.2.4}, denoted $R^{\text{BO}}_{\text{Opt}}$, can be written as
	\begin{align*}  \label{chap4eq:T-1.3}
		R^{\text{BO}}_{\text{Opt}} &= np[\frac{(1-p)}{p}t^{\text{BO}}_1+t^{\text{BO}}_2] \\
		&= np t^{\text{BO}}_2 (1+o_1(1)) \hspace{0.05cm},
		\tag{58}
	\end{align*}
where the term $o_1(1)$ is non-random, independent of index $i$ and tends to zero as $n \to \infty$.	Now note that under
    \hyperlink{assumption1}{Assumption 1}, $t_1/p \to 0$ as $n \to \infty$ where $t_1$ is the probability of type I error corresponding to $i^{\text{th}}$ decision rule of the form \eqref{eq:4.3.6} based on our class of one-group priors, as obtained in Theorem \ref{sparsecounttypeItun}. 
    Using this and Theorem \ref{sparsecounttypeIItun}, the
     asymptotic expression for the Bayes risk based on the decision rule, denoted by $R_{\text{OG}}$ can be written as
	\begin{equation*} \label{chap4eq:T-1.4}
		R_{\text{OG}} =npt_2(1+o_2(1)) \hspace{0.05cm},
		\tag{59}
	\end{equation*}
    by observing that $t_2$ is also bounded away from zero as in \eqref{chapbotypeIIlower} using the same argument.  Here also, the term $o_2(1)$ is non-random, independent of index $i$ and tends to zero as $n \to \infty$.
	Next, we take the ratio of \eqref{chap4eq:T-1.4} to \eqref{chap4eq:T-1.3}. The lower bound is obtained using the fact that $\frac{R_{\text{OG}}}{R^{\text{BO}}_{\text{Opt}}} \geq 1$. For the upper bound, observe that,
    the lower bound on 
    $t^{\text{BO}}_2$ is bounded away from zero
    and using Theorem \ref{sparsecounttypeIItun}, we have
    \begin{align*}
    \limsup_{n \to \infty}  \frac{R_{\text{OG}}}{R^{\text{BO}}_{\text{Opt}}} & \leq (\delta+1)^{\alpha} \mathbb{P} \bigg(Y \leq 2a+\alpha -\frac{2(a+\alpha)}{(\delta+2)}\bigg),
    \end{align*}
    where $Y \sim NB(\alpha,\frac{1}{\delta+1})$.
\end{proof}
\begin{proof}[Proof of Theorem \ref{sparsecountabosEB}]
	The asymptotic expression
	for the Bayes risk based on the decision rule \eqref{eq:4.3.8}, denoted by $R^{\text{EB}}_{\text{OG}}$, is of the form
	\begin{equation*} \label{chap4eq:T-2.1}
		R^{\text{EB}}_{\text{OG}}= \sum_{i=1}^{n} [(1-p)t^{\text{EB}}_{1i}+pt^{\text{EB}}_{2i}] = p\sum_{i=1}^{n} \bigg[\frac{1-p}{p}t^{\text{EB}}_{1i}+t^{\text{EB}}_{2i} \bigg] \hspace{0.05cm}.
		\tag{60}
	\end{equation*}
    Note that, we already proved in Theorem \ref{sparsecountabostun} that $t^{\text{BO}}_2=(\delta+\beta+1)^{-\alpha}$ is bounded away from $0$, for all sufficiently large $n$. Hence, using Theorem \ref{sparsecounttypeIIEB} along with \hyperlink{assumption1}{Assumption 1},
    the proof of our desired result is obtained if we establish that
	\begin{equation*} \label{chap4eq:T-2.2}
		\sum_{i=1}^{n} \frac{1-p}{p}t^{\text{EB}}_{1i}=o(n) \hspace{0.05cm} \text{as} \hspace{0.05cm} n \to \infty \hspace{0.05cm}.
		\tag{61}
	\end{equation*}
	To prove \eqref{chap4eq:T-2.2}, we use Theorem \ref{sparsecounttypeIEB} where, it is derived that the upper bound of $t^{\text{EB}}_{1i}$ is the same for each $i=1, 2,\cdots,n$. Further using the fact that $1-p \leq 1$, we obtain,
	\begin{equation*}  \label{chap4eq:T-2.3}
		\frac{1}{n}  \sum_{i=1}^{n} \frac{1-p}{p}t^{\text{EB}}_{1i} \leq \frac{2\alpha \beta}{a p}+\frac{\alpha \beta}{p}+\frac{1}{p} e^{-(2 \log 2-1)(1-(\beta+\delta+1)^{-\alpha})np(1+o(1))} \hspace{0.05cm}. 
		\tag{62}
	\end{equation*}
	By \hyperlink{assumption1}{Assumption 1}, the first two terms in the right-hand side of \eqref{chap4eq:T-2.3} go to $0$ as $n \to \infty$. For the third term, note that, for $p \propto n^{-\epsilon}, 0< \epsilon<1$, $np \to \infty$ as $n \to \infty$ and $\log(\frac{1}{p})=o(np)$. As a result, the third term too goes to $0$ as $n \to \infty$ and this proves \eqref{chap4eq:T-2.2}. This completes the proof of Theorem \ref{sparsecountabosEB}.
\end{proof}

\subsection{Distributions of different priors of the form \eqref{eq:4.1.4}}
\label{formpriors}
\begin{itemize}
    \item Three parameter beta normal mixtures: This class of priors can be written as
    \begin{align*} \label{countabostpbn}
        \pi_1(\lambda^2_i)= K (\lambda^2_i)^{-a_2-1} L(\lambda^2_i), \tag{63}
    \end{align*}
    where $L(\lambda^2_i)=(1+1/\lambda^2_i)^{-(a_1+a_2)}, K=\frac{\Gamma(a_1+a_2)}{\Gamma(a_1) \Gamma(a_2)}, a_1>0, a_2>0$. Hence, \eqref{countabostpbn} satisfies \eqref{eq:4.1.4} with $a=a_2$. 
    \item  Generalized double Pareto priors: This class of priors is of the form
    \begin{align*} \label{countabosgdp}
        \pi_1(\lambda^2_i)= K (\lambda^2_i)^{-a_1-1} L(\lambda^2_i), \tag{64}
    \end{align*}
 where $L(\lambda^2_i)= 2^{a_1-1} \int_{0}^{\infty} e^{-a_2 \sqrt{\frac{2z}{\lambda^2_i}}} e^{-z} z^{a_1} dz, a_1>0, a_2>0$.
 Therefore, \eqref{countabosgdp} can be expressed as \eqref{eq:4.1.4} with $a=a_1$.
\end{itemize}

\end{document}